\documentclass[12pt]{article}
\usepackage{graphicx,psfrag,epsfig}
\usepackage{amssymb,amsmath,amscd,amsthm}
\usepackage{graphicx,psfrag,epsfig}

\usepackage{graphicx}

\newtheorem{theorem}{Theorem}[section]

\newtheorem{lemma}[theorem]{Lemma}

\date{}

\begin{document}

\date{}
\title{Continuous Model for Homopolymers}
\author{M. Cranston\footnote{Dept of Mathematics, University of
California, Irvine, CA 92697, mcransto@math.uci.edu}, L.
Koralov\footnote{Dept of Mathematics, University of Maryland,
College Park, MD 20742, koralov@math.umd.edu}, S.
Molchanov\footnote{Dept of Mathematics, University of North
Carolina, Charlotte, NC 28223, smolchan@uncc.edu}, B.
Vainberg\footnote{Dept of Mathematics, University of North
Carolina, Charlotte, NC 28223, brvainbe@uncc.edu}} \maketitle

%\centerline{M. Cranston - Department of Mathematics, University of
%California} \centerline{Irvine, CA 92697, mcransto@math.uci.edu}
%
%\centerline{L. Koralov - Department of Mathematics, University of
%Maryland} \centerline{College Park, MD 20742,
%koralov@math.umd.edu}
%
%\centerline{S. Molchanov, B. Vainberg - Department of Mathematics,
%University of North Carolina} \centerline{Charlotte, NC 28223}

\begin{abstract}
We consider the model for the distribution of a long homopolymer
in a potential field. The typical shape of the polymer depends on
the temperature parameter. We show that at a critical value of the
temperature the transition occurs from a globular to an extended
phase. For various values of the temperature, including those at
or near the critical value, we consider the limiting behavior of
the polymer when its size tends to infinity.
\end{abstract}


\section{Introduction}

The goal of this paper is to analyze various critical phenomena
for a model of long homogeneous polymer chains in an attracting
potential field. The model exhibited here demonstrates a phase
transition from a densely packed globular phase at low
temperatures to an extended phase at higher temperatures. In the
latter phase, the thermal fluctuations overcome the attraction
between monomers and the chain takes on the shape of a $3d$ random
walk or Brownian motion with a typical scale $O(\sqrt{T})$ where
$T$ is the length of the polymer. A real life example of this
phenomenon is that of albumen (egg white). We describe a rough picture of this situation.
The physical reality is more complex as there are present several types
of protein with different critical points. However in a simplified version,
at room temperature the
albumen is in the globular state and as a result, it forms a
viscous, translucent liquid. However, at higher temperatures
(around $60-65^o$ C) there is a transition of the albumen to a
diffusive (extended) state resulting in an opaque semi-solid
material. While this transition may be reversible for an
individual polymer, in the aggregate, the polymer strands in the
diffusive state become interwoven and form chemical bonds with
each other and can not return to the globular state when the
temperature is decreased.

It is worthwhile recalling Gibbs' philosophy of phase transitions.
Start with a system of finite size $T$. The configuration space
$\Sigma_T=\{x(\cdot)\}$ denotes all possible states $x(\cdot)$ of
the system. The space $\Sigma_T$ is equipped with a reference
measure $\mathrm{P}_{0,T}$ which corresponds to infinite absolute
temperature (in our case, the inverse temperature $\beta=0$). The
configurations satisfy boundary conditions which reflect the
interaction of the finite system with its environment. This system
is endowed with a Hamiltonian $H_T$ giving the energy $H_T(x)$ of
the state $x$.  For $\beta>0,$ the Gibbs measure
$\mathrm{P}_{\beta,T}$ is given by the density
\begin{eqnarray}
\frac{d\mathrm{P}_{\beta,T}}{d\mathrm{P}_{0,T}}(x)=\frac{\exp(-\beta
H_T(x))}{Z_{\beta,T}},
\end{eqnarray}
where
\begin{eqnarray}
Z_{\beta,T}=\int_{\Sigma_T}\exp(-\beta H_T(x))dP_{0,T}.
\end{eqnarray}
When $T<\infty$, the measure $\mathrm{P}_{\beta,T}$ and the
thermodynamic quantities associated to $\mathrm{P}_{\beta,T}$ are
analytic functions of $\beta$.

Now let $T\to \infty.$ In typical situations, there is a critical
value $\beta_{cr}$ such that for  $\beta > \beta_{cr},$ there
exists a unique limiting measure $\mathrm{P}_\beta$ on $\Sigma,$
the space of infinite configurations, and this limiting measure is
independent of the boundary conditions on $\Sigma_T$. Moreover,
$\mathrm{P}_\beta$ and its relevant thermodynamic quantities are
still analytic functions of $ \beta$ for $\beta > \beta_{cr}.$
 One manifestation of the phase transition
 is the non-uniqueness for $\beta < \beta_{cr}$ of the limiting measure as
  $T \to \infty$ as it
   has dependence on the boundary conditions on $\Sigma_T$.
   Another is the non-analyticity of thermodynamic quantities
   associated to $\mathrm{P}_\beta$ as a function of  $\beta.$
   The mathematical characterization of the phase transition in terms of
   non-uniqueness of the limiting Gibbs measure traces its history to the
   works of Dobrushin \cite{Do} and Ruelle \cite{Ru}.

 Modern physical theories predict that near the critical
 point $\beta=\beta_{cr}$ the limiting Gibbs measure
 $\mathrm{P}_\beta$ must be invariant with respect to
 renormalizations of the system (self-similarity).
 This idea is related to the two-parametric scaling by
 Fisher \cite{Fi} for $\beta$ near $\beta_{cr}.$ Another important fact
 is that critical behavior as $\beta\to\beta_{cr}$ of the physical system demonstrates
 universality, that is the same behavior holds for a wide class of Hamiltonians.

 The most essential  part of the present paper is the detailed description of the
 polymer chain near the critical point and the establishment  of the physical
 ideas of universality and self-similarity for our particular model of homopolymers.

\section{Description of the Model and Results}

%In this paper we shall consider a model describing the
%distribution of a chain homopolymer in a potential field. Our
%approach is primarily based on investigating the asymptotic
%properties of solutions to parabolic equations.

A continuous function $x: [0,T] \rightarrow {\mathbb{R}}^d$, $x(0)
= 0$, will be thought of as a realization of the polymer. The
parameter $t \in [0,T]$ can be intuitively understood as the
length along the polymer (although the functions $x = x(t)$ are
not differentiable and the genuine notion of length can not be
defined).

We assume that for $\beta =0$,  the polymer is distributed
according to the Wiener measure $\mathrm{P}_{0,T}$ on $\Sigma_T =
C([0,T], {\mathbb{R}}^d)$. For an infinitely smooth compactly
supported potential $v \in C_0^\infty ( {\mathbb{R}}^d)$ and a
coupling constant $\beta \geq 0$, the polymer is distributed
according to the Gibbs measure $\mathrm{P}_{\beta,T}$, whose
density with respect to $\mathrm{P}_{0,T}$ is
\[
\frac{d \mathrm{P}_{\beta,T}}{d \mathrm{P}_{0,T}}(x) =
 \frac{\exp(\beta \int_0^T v(x(t)) d t)}{%
Z_{\beta,T}},~~x \in C([0,T], {\mathbb{R}}^d),
\]
In other words, the Hamiltonian $H_T$ is given by $H_T=-\int_0^T
v(x(s))d s$. The normalizing factor $Z_{\beta,T}$, called the
partition function, is given by
\begin{equation} \label{partf}
{Z_{\beta,T}} = \int_{C([0,T], {\mathbb{R}}^d)} \exp(\beta
\int_0^T v(x(t)) d t) d \mathrm{P}_{0,T} (x) = \mathrm{E}_{0,T}
e^{-\beta H_T}.
\end{equation}
 It
will be usually assumed that the potential is nonnegative and not
identically equal to zero. We shall be interested in the prevalent
behavior of the polymer with respect to the measure
$\mathrm{P}_{\beta,T}$ as $T \rightarrow \infty$.

 We shall see that
there are two qualitatively different cases corresponding to
different values of $\beta$. Namely, for all
sufficiently large values of $\beta$ there is a limiting distribution for $%
x(T)$ with respect to~$\mathrm{P}_{\beta,T}$. Moreover, for each
positive constant $s$ and each function $S(T)$ such that
$S(T)\rightarrow \infty $ and $ T - S(T) \rightarrow \infty$ as $T
\rightarrow \infty$, the family of processes $x(S(T)+ t)$, $t \in
[0, s]$, with respect to either measure  $ \mathrm{P}_{\beta,T}$
or $\mathrm{P}_{\beta,T}(\cdot|x(T)=0)$, converges to a Markov
process as $T \rightarrow \infty$. The generator of the limiting
Markov process and its invariant measure are written out
explicitly in Theorem~\ref{pro1}. Since $x(S(T))$ and $x(T)$
converge to limiting distributions and thus typically remain
bounded as $T \rightarrow \infty$, we shall say that the polymer
is in the globular state.

If $\beta > 0$ is sufficiently small and $d \geq 3$, then the
family of processes $x(t T)/\sqrt{T}$, $0 \leq t \leq 1$, defined
on $( C([0,T], {\mathbb{R}}^d), \mathrm{P}_{\beta,T})$, converges
to a Brownian motion on the interval $[0,1]$ (Theorem~\ref{nn1}).
In this case we shall say that the polymer is in the diffusive
state. Similarly, the family of processes $x(t T)/\sqrt{T}$, $0
\leq t \leq 1$, defined on $(( C([0,T], {\mathbb{R}}^d),
\mathrm{P}_{\beta,T}(\cdot|x(T)=0))$, converges to a Brownian
bridge on the interval $[0,1]$.

We shall see that there is a number $\beta_{cr}$ (called the critical value
of the coupling constant) such that the polymer is in the diffusive state
for $\beta < \beta_{cr}$ and in the globular state for $\beta > \beta_{cr}$.
The value of $\beta_{cr}$ and the behavior of the polymer when $\beta$ is
near $\beta_{cr}$ depend on the dimension $d$ and on the potential. In
particular, we shall see that $\beta_{cr} = 0$ for $d = 1,2$ and $\beta_{cr}
> 0$ for $d \geq 3$.

Of particular interest is the behavior of the polymer when $\beta
= \beta_{cr}$. In this case the appropriate scaling is the same as
in the diffusive case, that is we study the family of processes
$x(t T)/\sqrt{T}$, $0 \leq t \leq 1$. We shall find the limit of
this family as $T \rightarrow \infty$. It turns out to be a Markov
process with a non-Gaussian, spherically symmetric transition function
(Theorem~\ref{nn1cc}). The transition function of the limiting
Markov process will be written out explicitly.
%diffusion process
%with spherically symmetric generator (Theorem~\ref{nn1cc}). The
%coefficient at the diffusion term of the limiting process is
%constant, while the drift depends on both time and space and is
%directed towards the origin. The drift term becomes negligible
%when $t \uparrow 1$. The radial part of the limiting process is a
%diffusion process with reflection at the origin.

In order to determine whether the polymer is in the globular or diffusive
state for a given $\beta$, we shall look at the rate of growth of the
partition function ${Z_{\beta,T}}$. Namely, let
\[
\lambda_0(\beta) = \lim_{T \rightarrow \infty} \frac{\ln Z_{\beta,T}}{T}.
\]
It will be demonstrated that the limit exists and is equal to the
supremum of the spectrum of the operator $H_\beta =
\frac{1}{2}\Delta + \beta v: L^2({\mathbb{R}} ^{d})\rightarrow
L^2(\mathbb{R}^{d})$. The infimum of the set of $\beta$ for which
$\lambda_0(\beta) > 0$ is equal to $\beta_{cr}$. It will be seen
that $ \lambda_0(\beta_{cr}) = 0 $ is an eigenvalue of
$H_{\beta_{cr}}$ in dimensions $d \geq 5$, and corresponds to a
ground state of $H_{\beta_{cr}}$ in dimensions $d = 3,4$.

The paper is organized as follows.

In Section~\ref{timp} we consider finite $T$ and show that
$\{x(t), 0 \leq t \leq T\}$ is a time-inhomogeneous Markov process
with respect to the measures $\mathrm{P}_{\beta,T}$ and
$\mathrm{P}_{\beta,T}\left(\cdot| x(T)=0\right)$.

In Section~\ref{cvcc} we prove the existence of the critical value
of the coupling constant. In Section~\ref{sectt} we analyze the
properties of the resolvent of the operator $H_\beta$ which, in
particular, will be needed to study the asymptotic properties of
the partition function.

%As in other statistical physics models, of interest is the growth
%rate of $Z_{\beta, T}$ when $T \rightarrow \infty$ and
% $\beta$ is near or at $\beta_{cr}$.
In Section~\ref{secff} we shall examine the asymptotics of
$\lambda_0(\beta)$ when $ \beta \downarrow \beta_{cr}$
and show it has the following asymptotic behavior as $\beta\downarrow \beta_{cr},
$
\begin{eqnarray*}
\lambda_0(\beta)\sim\left\{\begin{array}{lll}c_3(\beta-\beta_{cr})^2,\,d=3,\\
c_4(\beta-\beta_{cr})/\ln(1/(\beta-\beta_{cr})),\,\,d=4,\\
c_d(\beta-\beta_{cr}),\,\,d\ge5.\end{array}\right.
\end{eqnarray*}
These asymptotics demonstrate universality in that they depend
only on dimension. The constants $c_d$, $d \geq 3$,  are not
universal however. In Section~\ref{secfi} we find the asymptotics,
as $T \rightarrow \infty$, of $Z_{\beta, T}$. In particular, when
$\beta>\beta_{cr}$, we shall find that $ Z_{\beta, T}\sim k_\beta
e^{\lambda_0(\beta)T}$ for some constant $k_\beta$, while for
$\beta<\beta_{cr},\,\,Z_{\beta, T}$ has a finite limit as
$T\to\infty.$ Finally, when $\beta=\beta_{cr},$ it turns out that
$Z_{\beta, T}\sim k_3 T^{1/2}$ for $d =3$,  $Z_{\beta, T}\sim k_4
{T}/{\ln T}$ for $d =4$, while $Z_{\beta, T}\sim k_d T$ for $d
\geq 5$. We also give asymptotics of the solutions to the
parabolic equation $\partial u / \partial t = H_\beta u$.
%In the case $\beta
%> \beta_{cr}$, the exponential rate of growth of the partition
%function and the solutions is equal to the principal eigenvalue
%$\lambda_0(\beta)$ studied in Section~\ref{secff}. For $\beta <
%\beta_{cr}$, $Z_{\beta, T}$ converges to a finite limit as $T
%\rightarrow \infty$. For $\beta = \beta_{cr}$, the growth rate of
%$Z_{\beta, T}$ depends on the dimension of the space. In
%particular, $Z_{\beta, T}$ grows as $T^{1/2}$ if $d = 3$.

In Sections~\ref{bp}, \ref{bb1x} and \ref{bb2},  we describe the
behavior of the polymer for $\beta > \beta_{cr}$, $\beta <
\beta_{cr}$ and $\beta = \beta_{cr}$, respectively, establishing
the convergence results mentioned above.

Some of the results presented  above have been obtained by
Cranston and Molchanov in \cite{CM} for the discrete model with
the potential concentrated at one point. The analysis was based on
explicit formulas for the solution of the parabolic equation with
such a potential. The current results demonstrate that the
behavior of the polymer is ``universal" with respect to the choice
of the potential. Another essential feature of this paper is the
detailed analysis of the behavior of the polymer when $\beta =
\beta_{cr}$. We refer the reader to the review of Lifschitz,
Grosberg and Khokhlov \cite{LGK} for a wealth of information and
ideas on polymer chains.

\section{Time-inhomogeneous Markov Property} \label{timp}

First we define $p_\beta$ as the fundamental solution of the heat
equation
\begin{eqnarray}
\begin{split}\label{pbeta}
\frac{\partial p_\beta}{\partial t}(t,y,x) =&  \frac{1}{2}\Delta_x
p_\beta(t,y,x) +\beta v(x) p_\beta(t,y,x),\\ p_\beta(0,y,x) =&
\delta(x-y).
\end{split}
\end{eqnarray}
In this section we shall prove that with respect to the measure
 $P_{\beta,T},$ the process $\{x(t), 0\le t\le T\}$ is a
 time-inhomogeneous Markov process.
Since we shall point out the link between non-uniqueness of Gibbs
measures and phase transitions it will be necessary to also
consider the transition mechanism for the process $\{x(t), 0 \leq
t \leq T \}$ under the conditional measure
$\mathrm{P}_{\beta,T}\left(\cdot|x(T)=0\right)$. Namely, we will
show that the free boundary condition corresponding to the measure
$\mathrm{P}_{\beta,T}$ and the pinned boundary condition
corresponding to the measure $\mathrm{P}_{\beta,T}\left(\cdot
|x(T)=0\right)$ lead to different Gibbs measures in the limit.

Let  $Z_{\beta,t}(x)=\mathrm{E}^x \exp(\beta \int_0^t v (x_s)d
s)$, where $\mathrm{E}^{x}$ is the expectation with respect to the
measure induced by the Brownian motion starting at $x$. Thus $
Z_{\beta,t}(0) = Z_{\beta,t}$, where $Z_{\beta,t}$ is the
partition function introduced in the previous section.
\begin{theorem}
The process  $\{x(t),0 \leq t \leq
T\}$ is a time-inhomogeneous Markov process with respect to the
measures $\mathrm{P}_{\beta,T}.$ Its transition density is given by
\begin{eqnarray}
q_\beta^T((s,y),(t,x))={p_\beta(t-s,y,x)Z_{\beta,T-t}(x)}(Z_{\beta,T-s}(y))^{-1},\label{qdens}.
\end{eqnarray}
The transition density $q_\beta^T((s,y),(t,x))$ solves the parabolic equation
\begin{eqnarray}
\frac{\partial}{\partial s}q^T_\beta((s,y),(t,x))+&\frac12
\Delta_y q^T_\beta((s,y),(t,x))+\nabla_y \ln Z_{\beta,T-s}(y)
\nabla_y q^T_\beta((s,y),(t,x))=0.
\end{eqnarray}
With respect to the conditional measure $\mathrm{P}_{\beta,T}(\cdot\,\,|x(T)=0),$
the process  $\{x(t),0 \leq t \leq
T\}$ is a time-inhomogeneous Markov process with transition density
\begin{eqnarray}
q_\beta^{(T,0)}((s,y),(t,x))=&p_\beta(t-s,y,x)p_\beta(T-t,x,0)(p_\beta(T-s,y,0))^{-1}.
\end{eqnarray}
\end{theorem}

While this
result is not used directly in later sections, it provides some
intuition on the nature of the limiting processes when we consider
the limit $T \rightarrow \infty$.

\begin{proof}
The Feynman-Kac formula gives that for $0<t\leq T$,
\begin{eqnarray}
\mathrm{P}_{\beta,T}(x(t)\in
dx)=\frac{p_\beta(t,0,x)\mathrm{E}^{x}\exp(\beta
\int_0^{T-t}v(x_s) d s)}{Z_{\beta,T}}d x.
\end{eqnarray}
 Similarly,
for $0 = t_0 < t_1< t_2 < ...< t_n \leq T$ and $x_0 = 0$,
\[
\mathrm{P}_{\beta,T}\left(x(t_1)\in dx_1,...,x(t_n)\in
dx_n\right)=
\]
\[
\frac{\prod_{i=0}^{n-1}
p_\beta(t_{i+1}-t_i,x_i,x_{i+1})\mathrm{E}^{x_n}\exp(\beta
\int_0^{T-t_n}v(x_s)d s)}{Z_{\beta,T}}dx_1dx_2...dx_n.
\]
So, if we  set for $0 \leq s < t \leq T$,
\[
q_\beta^T((s,y),(t,x))={p_\beta(t-s,y,x)Z_{\beta,T-t}(x)}(Z_{\beta,T-s}(y))^{-1},
\]
then
\[
\mathrm{P}_{\beta,T}\left(x(t_1)\in dx_1,...,x(t_n)\in
dx_n\right)=\prod_{i=0}^{n-1}q_{\beta}^T((t_i,x_i),(t_{i+1},x_{i+1})).
\]
Since $q_\beta^T((s,y),(t,x))>0$ and
\[
\int_{{{\mathbb{R}^d}}}q_\beta^T((s,y),(t,x))d x=1,
\]
 this means
that $\{x(t), 0 \leq t \leq T\}$ under the measure $P_{\beta,T}$
is a time-inhomogeneous Markov process with transition
probabilities $q^T.$ Turning the equation for $q_\beta^T$ around
and solving for $p_\beta$ yields
\[
p_\beta(t-s,y,x)=\frac{q_\beta^T((s,y),(t,x))Z_{\beta,T-s}(y)}{Z_{\beta,T-t}(x)}.
\]
Using the fact that
\[
\frac{\partial}{\partial
s}p_\beta(t-s,y,x)+\frac12\Delta_yp_\beta(t-s,y,x)+\beta
v(y)p_\beta(t-s,y,x)=0,
\]
we derive that $q_\beta^T$ satisfies the equation
\begin{eqnarray}
\begin{split}
\frac{\partial}{\partial s}q^T_\beta((s,y),(t,x))
\frac{Z_{\beta,T-s}(y)}{Z_{\beta,T-t}(x)}+&q^T_\beta((s,y),(t,x))\frac{\frac{\partial}{\partial
s}Z_{\beta,T-s}(y)}{Z_{\beta,T-t}(x)}\\ +\frac12 \Delta_y
q^T_\beta((s,y),(t,x)) \frac{Z_{\beta,T-s}(y)}{Z_{\beta,T-t}(x)}
+&\beta
v(y)q^T_\beta((s,y),(t,x))\frac{Z_{\beta,T-s}(y)}{Z_{\beta,T-t}(x)}\\
+\frac12\frac{q^T_\beta((s,y),(t,x))}{Z_{\beta,T-t}(x)}\Delta_y
Z_{\beta,T-s}(y)+&\nabla_y q^T_\beta((s,y),(t,x))\frac{\nabla_y
Z_{\beta,T-s}(y)}{Z_{\beta,T-t}(x)}\\ =&0.
\end{split}
\end{eqnarray}
Simplifying this leads to the following parabolic equation for
$q_\beta^T$,
\begin{eqnarray}
\frac{\partial}{\partial s}q^T_\beta((s,y),(t,x))+&\frac12
\Delta_y q^T_\beta((s,y),(t,x))+\nabla_y \ln Z_{\beta,T-s}(y)
\nabla_y q^T_\beta((s,y),(t,x))=0.
\end{eqnarray}

Next we consider the pinned case,  for $0 = t_0 < t_1<...< t_n <t_{n+1}=T$ and
 $x_0=x_{n+1}=0.$ Then,
\begin{eqnarray}
\begin{split}
\mathrm{P}_{\beta,T}\left(x(t_1)\in dx_1,...,x(t_n)\in
dx_n|x(T)=0\right)=&\frac{\mathrm{P}_{\beta,T}\left(x(t_1)\in
dx_1,...,x(t_n)\in
dx_n,x(T)=0\right)}{\mathrm{P}_{\beta,T}\left(x(T)=0\right)}\\
=&\frac{\prod_{i=0}^{n-1}
p_\beta(t_{i+1}-t_i,x_i,x_{i+1})}{p_\beta(T,0,0)}dx_1...dx_n.\\
\end{split}
\end{eqnarray}
Now  set for $0 \leq s<t \leq T$,
\begin{eqnarray}
q_\beta^{(T,0)}((s,y),(t,x))=&p_\beta(t-s,y,x)p_\beta(T-t,x,0)(p_\beta(T-s,y,0))^{-1}.
\end{eqnarray}
Then
\begin{eqnarray}
\mathrm{P}_{\beta,T}\left(x(t_1)\in dx_1,...,x(t_n)\in
dx_n|x(T)=0\right)=\prod_{i=0}^{n-1}q_{\beta}^{(T,0)}((t_i,x_i),(t_{i+1},x_{i+1})).
\end{eqnarray}
Since $q_\beta^{(T,0)}((s,y),(t,x))>0$ and
\[
\int_{{{\mathbb{R}^d}}}q_\beta^T((s,y),(t,x))d x=1,
\]
 this means
that $\{x(t), 0 \leq t \leq T\}$ under the conditional measure
$\mathrm{P}_{\beta,T}\left(\cdot | x(T)=0\right)$ is a
time-inhomogeneous Markov process with transition densities
$q_\beta^{(T,0)}$.
\end{proof}
We shall see below in that in the globular phase
$\beta>\beta_{cr}$ the drift term $\nabla_x \ln Z_{\beta,T-s}(x)$
has a non-trivial limit as $T\to\infty$. This means that for
$\beta>\beta_{cr},$ the Gibbs measure corresponds to a stationary
Markov process in the $T\to\infty$ limit. On the other hand, this
limit will vanish for $\beta<\beta_{cr}.$ This explains the nature
of the diffusive state for high temperature.

\section{Critical Value of the Coupling Constant}

\label{cvcc} Let
\[
H_{\beta }=\frac{1}{2}\Delta +\beta
v:L^{2}({\mathbb{R}}^{d})\rightarrow L^{2}({\mathbb{R}}^{d}),
\text{ \ \ \ }v=v(x)\in C_{0}^{\infty }(\mathbb{R}^{d}),~~~\beta
\geq 0.
\]
%Let $\Omega$ be an arbitrary closed ball containing the support of
%$v$.
 We shall always assume that $v(x)$ is non-negative and compactly supported,
 although
many results do not require these restrictions or can be modified
to be valid without these restrictions. We shall also assume that
$v$ is not identically equal to zero. It is well-known that the
spectrum of $H_{\beta }$ consists of the absolutely continuous part $%
(-\infty ,0]$ and at most a finite number of non-negative eigenvalues:
\[
\sigma (H_{\beta })=(-\infty ,0]\cup \{\lambda _{j}\},\text{ \ \ }0\leq
j\leq N,\text{ \ \ }\lambda _{j}=\lambda _{j}(\beta )\geq 0.
\]
We enumerate the eigenvalues in the decreasing order. Thus, if
$\{\lambda _{j}\}\neq \emptyset $, then $\lambda _{0}=\max \lambda
_{j}$.
\begin{lemma}
\label{lmonot} There exists $\beta _{cr} \geq 0$ (which will be
called the critical value of $\beta $) such that $\sup \sigma
(H_{\beta })=0$ for $\beta \leq \beta _{cr}$ and $\sup \sigma
(H_{\beta })=\lambda _{0}(\beta)>0$ for $\beta >\beta _{cr}$. For
$\beta
>\beta _{cr}$ the eigenvalue $\lambda _{0}(\beta )$ is a strictly
increasing and continuous function of $\beta $. Moreover,
$\lim_{\beta \downarrow \beta _{cr}}\lambda (\beta )=0$ and
$\lim_{\beta \uparrow \infty
}\lambda (\beta )=\infty $. %depends
%monotonically on $\beta ,$ i.e. if $\beta _{2}>\beta _{1}$ and
%$\lambda _{0}(\beta _{1})$ exists then $\lambda _{0}(\beta _{2})$
%exists and $\lambda _{0}(\beta _{2})>\lambda _{0}(\beta _{1}).$
\end{lemma}
\proof
 The form $(H_{\beta }\psi ,\psi )$ is positive on a
function $\psi $ supported on ${\rm supp} (v) $ if $\beta $ is
large enough. Thus $\sup \sigma (H_{\beta })>0$ for sufficiently
large $\beta$. On the other hand, $\sigma (H_{\beta })=(-\infty
,0]$ when $\beta =0$. Let $\beta_{cr} = \sup\{ \beta: \sup \sigma
(H_{\beta })=0 \}$. It is clear that
 $\sup \sigma (H_{\beta })=0$ for
$\beta < \beta_{cr}$ since the operator $H_{\beta }$ depends
monotonically on $\beta$.

 Other statements  easily follow
from the fact that for each $\psi$ the form $(H_{\beta }\psi ,\psi
)$ depends continuously and monotonically on $\beta$. \qed
\\

\noindent {\bf Remark.} As will be shown below, $\beta_{cr} = 0$
for $d =1,2$, and $\beta_{cr} \geq 0$ for $d \geq 3$. Thus we do
not talk about phase transition for $d =1,2$ since we do not
consider negative values of $\beta$.

For $d \geq 3$, by the Cwikel-Lieb-Rozenblum estimate \cite{RS},
\[
\sharp \{ \lambda_i(\beta) \geq 0 \} \leq c_d \beta^{d/2} \int_{
\mathbb{R}^d } | v(x)|^{d/2} d x.
\]
This implies that there are no eigenvalues for sufficiently small
values of $\beta$ if $d \geq 3$, that is $\beta_{cr} > 0$. It is
also well-know (see \cite{RS}) that $\sup \sigma (H_{\beta }) > 0$
for $d =1,2$ if $\beta > 0$, $v \geq 0$ and $v$ is not identically
zero. These statements will also be proved below without referring
to the Cwikel-Lieb-Rozenblum estimate.

\section{Analytic Properties of the Resolvent} \label{sectt}

The resolvent  of the operator $H_{\beta }$ will be considered in
the spaces of square-integrable and continuous functions.  The
resolvent $R_{\beta }(\lambda )=(H_{\beta }-\lambda )^{-1}:
L^{2}({\mathbb{R}}^{d})\rightarrow L^{2}({\mathbb{R}}^{d})$
 is a meromorphic operator valued function on $\mathbb{C}^{\prime }
 =\mathbb{C}\backslash (-\infty ,0].$ Denote the kernel
of $R_{\beta }(\lambda )$ by $R_{\beta }(\lambda ,x,y).$ If $\beta
=0,$ the kernel depends on the difference $x-y$ and will be
denoted by $R_{0}(\lambda ,x-y).$
The kernel $R_{0}(\lambda ,x)$ can be expressed through the Hankel function $%
H_{\nu }^{(1)}$:
\begin{equation}  \label{ha1}
R_{0}(1,x)=c|x|^{1-\frac{d}{2}}H_{\frac{d}{2}-1}^{(1)}(i\sqrt2|x|),
\end{equation}
and
\begin{equation}  \label{ha2}
R_{0}(\lambda ,x)=
ck^{d-2}(k|x|)^{1-\frac{d}{2}}H_{\frac{d}{2}-1}^{(1)}(i\sqrt2k|x|),~~
k=\sqrt{\lambda},~~\mathrm{Re}k > 0.
\end{equation}
In particular,
\[
R_{0}(\lambda ,x)=\frac{e^{-\sqrt2k|x|}}{-\sqrt2k},~~~d=1;~~~R_{0}(\lambda ,x)=
\frac{e^{-\sqrt2k|x|}}{-2\pi|x|},~~~d=3.
\]
We shall say that $f \in L^2_{\exp}( \mathbb{R}^d)$ if $f$ is
measurable and
\[
||f||_{ L^2_{\exp}( \mathbb{R}^d)} = (\int_{ \mathbb{R}^d} f^2(x)
e^{ |x|^2} d x)^{\frac{1}{2}} < \infty.
\]
Similarly, we shall say that  $f \in C_{\exp}( \mathbb{R}^d)$ if
$f$ is continuous and
\[
||f||_{C_{\exp}( \mathbb{R}^d)} = \sup_{x \in \mathbb{R}^d}(
|f(x)| e^{{|x|^2}}) < \infty.
\]

Note that $R_{0}(\lambda )$, $\lambda \in \mathbb{C}'$, is a
bounded operator not only in $L^2( \mathbb{R}^d)$ but also from
$C_{\exp}( \mathbb{R}^d)$ to $C( \mathbb{R}^d)$, where $C(
\mathbb{R}^d)$ is the space of bounded continuous functions on
$\mathbb{R}^d$. Denote
\begin{equation}
A(\lambda )= v(x)R_{0}(\lambda ):~ L^2_{\exp}( \mathbb{R}^d)
\rightarrow L^2_{\exp}( \mathbb{R}^d)~~ ({\rm and}~~ C_{\exp}(
\mathbb{R}^d) \rightarrow C_{\exp}( \mathbb{R}^d)). \label{al1}
\end{equation}

 The well-known properties of the
Hankel functions together with (\ref {ha1}) and (\ref{ha2}) imply
the following lemma (see \cite{Va} for a similar statement for
general elliptic operators).
\begin{lemma}
\label{ker10} Consider the operator $A(\lambda )$ in the spaces
$L^2_{\exp}(\mathbb{R}^d)$ and $C_{\exp}(\mathbb{R}^d)$.

(1) The operator $A(\lambda )$  is analytic in $\lambda \in
\mathbb{C}'$. It   admits an analytic extension as an entire
function of $\sqrt{\lambda }$ if $d$ is odd, except $d=1$, when it
has a pole (with respect to $\sqrt{\lambda }$) at the origin. The
operator $A(\lambda )$ has the form $A(\lambda )=A_{1}(\lambda
)+\ln \lambda A_{2}(\lambda )$ if $d$ is even, where $A_{1}$ and
$A_{2}$ are entire functions.

(2) $A_{2}(0)=0$ if $d \geq 4$ ($d$ is even), and therefore
$A(0)=\lim_{\lambda \rightarrow 0, \lambda \in
\mathbb{C}'}A(\lambda)$ exists and is a bounded operator for all
$d\geq 3$.

(3) The operator $A(\lambda )$ is compact for all $\lambda \in
\mathbb{C}' \cup \{0\}$ ($\lambda \neq 0$ if $d=1$ or $2$).

 (4) For each $\varepsilon
> 0$, we have $||A(\lambda)||=O(1/|\lambda|)$ as $\lambda \rightarrow
\infty$, $|{\rm arg} \lambda| \leq \pi - \varepsilon$.

 (5
) The operator
$A(\lambda )$ has the following asymptotic behavior as $\lambda
\rightarrow 0$, $\lambda \in \mathbb{C}' $:
\[
A(\lambda )=-{v P_{1}}/{\sqrt{\lambda }}+O(1),~~~ d=1,
\]
\[
 A(\lambda )=-{vP_{2}}{\ln ({1}/{\lambda}) }+O(1),~~~d=2,
\]
\[
A(\lambda )=-v({P_{3}}+Q_{3}\sqrt{\lambda })+O(|\lambda |),~~~d=3,
\]
\[
A(\lambda )=-v({P_{4}}+Q_{4}\lambda \ln ({1}/{\lambda}))
+O(|\lambda |),~~~d=4,
\]
\[
A(\lambda )=-v({P_{d}}+Q_{d}\lambda) +O(|\lambda |^{3/2}),~~~d\geq
5,
\]
where the operators $P_{d}$, $d\geq 1$,  $Q_{d}$, $d \geq 3$, have
the following kernels:
\[
P_{1}(x,y)=\frac{1}{\sqrt2 },~~~P_{2}(x,y)=\frac{1}{\pi },
\]
\[
P_{3}(x,y)=\frac{1}{2\pi |x-y|},~~~Q_{3}(x,y)=-\frac{1}{\sqrt2\pi },
\]
\[
P_{4}(x,y)=\frac{1}{\pi^2
|x-y|^2},~~~Q_{4}(x,y)=-\frac{1}{2\pi^2},
\]
\[
P_{d}(x,y)=\frac{a_{d}}{|x-y|^{d-2}},~~~Q_{d}(x,y)=\frac{-a_{d}}{%
(d-4)|x-y|^{d-4}},~~a_d > 0,~d \geq 5.
\]
\end{lemma}
\proof Let $d$ be odd.
% Consider the operator $A(\lambda )$  in
% $C_{\exp}(\mathbb{R}^d)$ (the case of $L^2_{\exp}(\mathbb{R}^d)$
% is simpler).
From (\ref{ha1}), (\ref{ha2}) and (\ref{al1}) it follows that the
kernel $A(\lambda ,x,y)=v(x)R_0(\lambda , x-y)$ of the operator
${A}(\lambda )$ is an entire function of $k= \sqrt \lambda$ if $d
\geq 3$ (but has a pole at $k=0$ if $d=1$). The kernel has a weak
singularity at $x=y$ and an exponential estimate at infinity. To
be more exact,
\begin{equation} \label{estt}
|A(k^2 ,x,y)|+|\frac{\partial A(k^2 ,x,y)}{ \partial k}| \leq
C(d,k)|v(x)|e^{|k(x-y)|}(|x-y|+|x-y|^{-(d-1)}),
\end{equation}
were $C(d,k)$ has a singularity at $k=0$ if $d=1$. Since
\[
||A(k^2)||_{C_{\exp}(\mathbb{R}^d)} \leq \sup_{x \in \mathbb{R}^d}
\int e^{|x|^2-|y|^2}|A(k^2 ,x,y)|d y,
\]
\[
 ||\frac {d}{d k} A(k^2)||_{C_{\exp}(\mathbb{R}^d)} \leq
\sup_{x \in \mathbb{R}^d} \int e^{|x|^2-|y|^2}| \frac{\partial
A(k^2 ,x,y)}{
\partial k}|d y,
\]
the estimate (\ref{estt}) immediately leads to the analyticity in
$k= \sqrt \lambda$ of the operator $A(\lambda)$ in the space
$C_{\exp}(\mathbb{R}^d)$. In order to get the same result in the
space $L^2_{\exp}(\mathbb{R}^d)$, we represent $A(\lambda)$ in the
form $B_1+B_2$ were the kernel $B_1(\lambda,x,y)$ of the operator
$B_1$ is equal to $\chi(x-y)A(\lambda, x, y)$. Here $\chi$ is the
indicator function of the unit ball. Since
\begin{equation}
|\chi(x) R_0(k^2, x)| + |\frac {d}{d k} \chi(x) R_0(k^2, x)| \in
L^1(\mathbb{R}^d),
\end{equation}
the convolution with $\chi(x) R_0(k^2, x)$ is an analytic in $k$
operator in the space $L^2(\mathbb{R}^d)$. Then $B_1$ (which is
the convolution followed by multiplication by $v(x)$) is an
analytic operator in the space $L^2_{\exp}(\mathbb{R}^d)$. The
product of the kernel of the operator $B_2$ and $e^{|x|^2-|y|^2}$
is square integrable in $(x,y)$. The same is true for the
derivative in $k$ of the kernel of $B_2$ multiplied by
$e^{|x|^2-|y|^2}$. Thus $B_2$ is also analytic in $k$. This
completes the proof of the analyticity of $A(\lambda)$ when $d$ is
odd. The case of even $d$ is similar. One needs only to take into
account that $R_0(\lambda ,x)$ has a logarithmic branching  point
at $\lambda =0$ in this case. The second statement of the lemma
follows immediately from (\ref{ha1}), (\ref{ha2}) and (\ref{al1}).

To prove the compactness of $A(\lambda)$,
we note that the estimate (\ref{estt})
 is valid not only for $A(k^2 ,x,y)$ and $\partial A(k^2 ,x,y) / \partial k$,
  but also for $\nabla_x A(k^2 ,x,y)$.
  Thus the arguments above lead to
  the boundedness of the operators $\frac {\partial}{\partial x_i}A(\lambda)$
(the composition of $A(\lambda)$ with the differentiation). Since
the supports of functions $A(\lambda)f$ belong to the support of
$v$, the standard Sobolev embedding theorems imply the compactness
of the operator $A(\lambda)$ in both the spaces
$L^2_{\exp}(\mathbb{R}^d)$ and $C_{\exp}(\mathbb{R}^d)$.

In order to prove the fourth statement of the lemma, we observe
that the $L^2( \mathbb{R}^d)$ norm of the resolvent $R_0(\lambda)$
does not exceed $1/| \rm {Im} \lambda|$ (the inverse distance from
the spectrum). Since $A(\lambda)$ is obtained from $R_0(\lambda)$
after multiplying it by a bounded function with compact support,
the $L^2_{\exp}( \mathbb{R}^d)$ norm of $A(\lambda)$ does not
exceed $c/| \rm {Im} \lambda|$, where $c$ is a positive constant
which depends on $v$.
 The norm of $A(\lambda)$ in the space $C_{\exp}(
\mathbb{R}^d)$ can be estimated by $\sup_{x \in \mathbb{R}^d}
|v(x) e^{|x|^2}| \int_{ \mathbb{R}^d} |{R}_{0}(\lambda ,x)| d x$,
which is of order $O(1/|\lambda|)$ as
 $\lambda \rightarrow \infty$, $|{\rm arg}
\lambda| \leq \pi - \varepsilon$, due to (\ref{ha2}).

 The
remaining statements  also easily follow from (\ref{ha1})
and~(\ref{ha2}). \qed

 Note that for $d \geq 3$, there
exists the limit
\[
R_{0}(0 ,x-y) := \lim_{\lambda \rightarrow 0, \lambda \in
\mathbb{C}' }R_{0}(\lambda ,x-y)=-a_{d}|x-y|^{2-d},
\]
which is a fundamental solution of the operator
$\frac{1}{2}\Delta$. The operator with this kernel will be denoted
by $R_0(0)$. While $R_0(\lambda)$, $\lambda \in \mathbb{C}'$, acts
in $L^2(\mathbb{R}^d)$ and $C( \mathbb{R}^d)$, the operator
$R_0(0)$ only maps $C_{\exp}( \mathbb{R}^d)$ to $C( \mathbb{R}^d)$
if $d < 5$. The following lemma follows from formulas (\ref{ha1})
and~(\ref{ha2}) similarly to Lemma~\ref{ker10}.
\begin{lemma} \label{conti1} For $d \geq 3$, the operator
$R_0(\lambda)$ considered as an operator from $C_{\exp}(
\mathbb{R}^d)$ to $C( \mathbb{R}^d)$ is analytic in $\lambda \in
\mathbb{C}'$.   It is uniformly bounded in $ \mathbb{C}'$. For
each $\varepsilon > 0$, it is of order $O(1/|\lambda|)$ as
$\lambda \rightarrow \infty$, $|{\rm arg} \lambda| \leq \pi -
\varepsilon$. It has the following asymptotic behavior as $\lambda
\rightarrow 0$, $\lambda \in \mathbb{C}'$:
\[
{R}_{0}(\lambda )={R}_{0}(0)+O(\sqrt{|\lambda|}),~~~d=3,
\]
\[
 {R}_{0}(\lambda
)={R}_{0}(0) +O(|\lambda \ln \lambda |),~~~d=4,
\]
\[
{R}_{0}(\lambda )={R}_{0}(0) +O(|\lambda |),~~~d\geq 5.
\]
\end{lemma}

The following lemma is simply a resolvent identity. It plays an
important role in our future analysis.

\begin{lemma}
For $\lambda \in \mathbb{C}^{\prime}$, we have the following
relation between the meromorphic operator-valued functions
\label{lrid}
\begin{equation}  \label{aa1}
R_{\beta}(\lambda)=R_0(\lambda) - R_0(\lambda)(I+\beta
v(x)R_0(\lambda))^{-1}[\beta v(x) R_0(\lambda)]
\end{equation}
\end{lemma}

\noindent {\bf Remark.} Note that (\ref{aa1}) can be written as
\begin{equation}
R_{\beta }(\lambda )=R_{0}(\lambda )-R_{0}(\lambda )(I+\beta
A(\lambda ))^{-1}[\beta v(x)R_{0}(\lambda )].  \label{aa2}
\end{equation}
From here it also follows that
\begin{equation}
R_{\beta }(\lambda )=R_{0}(\lambda )(I+\beta A(\lambda ))^{-1},  \label{aa3}
\end{equation}
which should be understood as an identity between meromorphic in
$\lambda $ operators acting from $L^2_{\exp}(\mathbb{R}^d)$ to
$L^2(\mathbb{R}^d)$ and from $C_{\exp}( \mathbb{R}^d)$ to $C(
\mathbb{R}^d)$. In the lattice case considered in~\cite{CM}, the
operator $A(\lambda)$ has rank one and
\[
R_\beta(\lambda,x,y)=R_0(\lambda,x,y)/(1-\beta I(\lambda)),\]
where $I(\lambda)$ is an analytic function of $\sqrt{\lambda}$
related to $A(\lambda).$ This exact formula is the key to all the
results in \cite{CM}.
\\

The kernels of the operators $I+\beta A(\lambda )$ (both in spaces
$L^2_{\exp}(\mathbb{R}^d)$ and $C_{\exp}( \mathbb{R}^d)$) are
described by the following lemma.
\begin{lemma}
\label{re1} (1) The operator-valued function $(I + \beta
A(\lambda))^{-1}$ is meromorphic in $\mathbb{C}^{\prime}$. It has
a pole at $\lambda \in \mathbb{C}^{\prime}$ if and only if
$\lambda$ is an eigenvalue of $H_\beta$. These poles are of the
first order.

(2) Let $\lambda_i(\beta)$ be a positive eigenvalue of $H_\beta$. There is a
one-to-one correspondence between the kernel of the operator $I + \beta
A(\lambda_i)$ and the eigenspace of the operator $H_\beta$ corresponding to
the eigenvalue $\lambda_i$. Namely, if $(I + \beta A(\lambda_i))h = 0$, then
$\psi =- R_0(\lambda_i) h$ is an eigenfunction of $H_\beta$ and $h = \beta v
\psi$.

(3) If $d \geq 3$, there is a one-to-one correspondence between
the kernel of the operator $I+\beta A(0)$ and solution space of
the problem
\begin{equation}
H_{\beta }(\psi )=\frac{1}{2}\Delta \psi +\beta v(x)\psi
=0,~~~\psi (x) = O(|x|^{2-d}),~~\frac{\partial{\psi}
}{\partial{r}}(x) = O(|x|^{1-d})~~as~~r =|x| \rightarrow \infty.
\label{grs}
\end{equation}
Namely, if $(I + \beta A(0))h = 0$ for $h \in L^2_{\exp}(
\mathbb{R}^d)$, then $h \in C_{\exp}( \mathbb{R}^d)$, $\psi
=-R_{0}(0)h$ is a solution of (\ref{grs}) and $h = \beta v \psi$.
\end{lemma}

\noindent {\bf Remark.} The relations (\ref{grs}) are an analogue
of the eigenvalue problem for zero eigenvalue and the
eigenfunction $\psi $ which does not necessarily belong to
$L^{2}(\mathbb{R}^d)$ (see Lemma~\ref{tgrs} below). We shall call
a non-zero solution of (\ref{grs}) a ground state.
% if $\beta = \beta_{cr}$ (and therefore
%$\lambda = 0$ is the edge of the spectrum of $H_\beta$).

 \proof
The operator $A(\lambda)$, $\lambda \in \mathbb{C}'$, is analytic,
compact, and tends to zero as $\lambda \rightarrow +\infty$ by
Lemma~\ref{ker10}. Therefore $(I+\beta A(\lambda ))^{-1}$ is
meromorphic  by the Analytic Fredholm Theorem.

 If $\lambda \in \mathbb{C}^{\prime }$ is a pole
of $(I+\beta A(\lambda ))^{-1}$, then it is also a pole of the
same order of $R_{\beta }(\lambda )$ as follows from (\ref{aa3})
since the kernel of $R_{0}(\lambda )$ is trivial. Therefore the
pole is simple and coincides with one of the eigenvalues $\lambda
_{i}$. Note that $\lambda $ is a pole of $(I+\beta A(\lambda
))^{-1}$ if and only if the kernel of $I+\beta A(\lambda )$ is
non-trivial. Let $h\in L^2_{\exp}(\mathbb{R}^d)$ be such that
$||h||_{L^{2}_{\exp}(\mathbb{R}^d)}\neq 0$ and $(I+\beta
vR_{0}(\lambda ))h=0$. Then $\psi :=-R_{0}(\lambda
)h\in L^{2}(\mathbb{R}^{d})$ and $(\frac{1}{2}\Delta -\lambda +\beta v)\psi =0$, that is $%
\psi $ is an eigenfunction of $H_{\beta }$.

Conversely, let $\psi \in L^{2}(\mathbb{R}^{d})$ be an
eigenfunction corresponding to an eigenvalue $\lambda_{i}$, that
is
\begin{equation}
(\frac{1}{2}\Delta -\lambda _{i})\psi +\beta v\psi =0.
\label{bb1}
\end{equation}
Denote $h=\beta v\psi $. Then  $(\frac{1}{2}\Delta -\lambda
_{i})\psi =-h$. Thus $\psi =-R_{0}(\lambda _{i})h$ and (\ref{bb1})
implies that $h$ satisfies $(I+\beta vR_{0}(\lambda _{i}))h=0$.
Note that $h \in C^\infty( \mathbb{R}^d)$, $h$ vanishes outside
${\rm supp} (v)$, and therefore belongs to the kernel of $I+\beta
A(\lambda _{i})$. This completes the proof of the first two
statements.

Similar arguments can be used to prove the last statement. If
$h\in L^2_{\exp}(\mathbb{R}^d)$ is such that
$||h||_{L^{2}_{\exp}(\mathbb{R}^d)}\neq 0$ and $(I+\beta
A(0))h=0$, then $h$ has compact support and  the integral operator
$R_0(0)$ can be applied to $h$. It is clear that $\psi
:=-R_{0}(0)h$ satisfies (\ref{grs}) and, since $h$ has compact
support, $h \in C_{\exp}( \mathbb{R}^d)$.

In order to prove that any solution of (\ref{grs}) corresponds to
an eigenvector of $I + \beta A(0)$, one only needs to show that
the solution $\psi $ of the problem (\ref{grs}) can be represented
in the form $\psi =-R_{0}(0)h$ with $h=\beta v\psi .$ The latter
follows from the Green formula
\[
\psi(x) =-(R_{0}(0)(\beta v\psi))(x)+\int_{|y|=a}[R_{0}(0,x-y)\psi'_{r}(y)-\frac{\partial }{%
\partial r}R_{0}(0,x-y)\psi (y)]ds,\text{ \ \ \ }|x|<a,
\]
after passing to the limit as $a\rightarrow \infty .$ \qed

Lemma~\ref{re1} can be improved for $\lambda = \lambda_0(\beta)$.
Due to the monotonicity and continuity of $\lambda =
\lambda_0(\beta)$ for $\beta > \beta_{cr}$, we can define the
inverse function
\begin{equation} \label{bla}
\beta = \beta(\lambda): [0, \infty) \rightarrow
[\beta_{cr},\infty).
\end{equation}

We shall prove that the operator $-A(\lambda)$, $\lambda
>0$, has a non-negative kernel and has a positive simple
eigenvalue such that all the other eigenvalues are smaller in
absolute value. Such an eigenvalue is called the principal
eigenvalue.

\begin{lemma} \label{ll44} The operator $-A(\lambda)$, $\lambda > 0$,
(in the spaces $L^2_{\exp}(\mathbb{R}^d)$ and $C_{\exp}(
\mathbb{R}^d)$) has the principal eigenvalue. This eigenvalue is
equal to $1/\beta(\lambda)$ and the corresponding eigenfunction
can be taken to be positive in the interior of ${\rm supp} (v)$
and equal to zero outside of ${\rm supp} (v)$. If $d \geq 3$, then
the same is true for the operator $-A(0)$ (in particular,
$\beta_{cr} > 0$).
\end{lemma}

\noindent {\bf Remark 1.}  Let $d \geq 3$.  Lemmas~\ref{re1} and
\ref{ll44} imply that the ground state of the operator $H_\beta$
for $\beta = \beta_{cr}$ (defined by (\ref{grs})) is defined
uniquely up to a multiplicative constant and corresponds to the
principal eigenvalue of $A(0)$. The ground state (with $\lambda =
0$) does not exist if $\beta < \beta_{cr}$.
\\

\noindent {\bf Remark 2.} Let $d \geq 3$.
 From Lemma~\ref{ker10} it follows that
\[
\lim_{\lambda \rightarrow 0, \lambda \in \mathbb{C}'} A(\lambda) =
A(0).
\]
Therefore for all  $\lambda \in \mathbb{C}'$ with $|\lambda|$
sufficiently small, the operator $-A(\lambda)$ has a simple
eigenvalue whose real part is larger than the absolute values  of
the other eigenvalues. We shall denote this eigenvalue by
$1/\beta(\lambda)$, thus extending the domain of the function
$\beta(\lambda)$ (see (\ref{bla})) from $[0,\infty)$ to $
[0,\infty) \cup (U \cap \mathbb{C}')$, where $U$ is a sufficiently
small neighborhood of zero.
\\

\noindent {\it Proof of Lemma~\ref{ll44}.} By Lemma~\ref{re1} it
is sufficient to consider the case of $L^2_{\exp}(\mathbb{R}^d)$.
The maximum principle for the operator $(\frac{1}{2}\Delta
-\lambda )$, $\lambda
>0$, implies that the kernel of the operator $R_{0}(\lambda )$,
$\lambda >0$, is negative. Thus, by (\ref{al1}), for all $y$ the
kernel of $-A(\lambda )$ is positive when $x$ is in the interior
of  ${\rm supp} (v) $ and zero otherwise. Thus $-A(\lambda )$,
$\lambda
>0$, has the principal eigenvalue (see \cite{Kr}).
On the other hand, by Lemma~\ref{re1}, $1/\beta(\lambda)$ is a positive eigenvalue of $%
-A(\lambda)$. Note that this is the largest positive eigenvalue of $%
-A(\lambda)$. Indeed, if $\mu = 1/\beta^{\prime}>
1/\beta(\lambda)$ is an eigenvalue of $-A(\lambda)$, then
$\lambda$ is one of the eigenvalues $\lambda_i$ of
$H_{\beta^{\prime}}$ by
Lemma~\ref{re1}. Therefore, $%
\lambda_i(\beta^{\prime}) = \lambda_0(\beta)$ for $\beta^{\prime}<
\beta$. This contradicts the monotonicity of $\lambda_0(\beta)$.
Hence the statement of the lemma concerning the case $\lambda
> 0$ holds.

For $d \geq 3$, the kernel of $-A(0)$ is equal to $v P_d$ and has
the same properties as the kernel of $-A(\lambda)$, $\lambda > 0$.
Thus $-A(0)$ has the principal eigenvalue. Since $A(\lambda)
\rightarrow A(0)$ as $\lambda \downarrow 0$, the principal
eigenvalue $1/\beta(\lambda)$ converges to the principal
eigenvalue $\mu < \infty$ of $-A(0)$. On the other hand,
$\beta(\lambda)$ is a continuous function, and therefore $\mu =
1/\beta_{cr}$, which proves the statement concerning the case
$\lambda = 0$. \qed
\\

The relationship between ground states and eigenfunctions of
$H_\beta$ is explained by the following lemma.

\begin{lemma}
\label{tgrs} Let $\beta =\beta _{cr}$.  If $d = 3$ or $d =4$, then
$H_{\beta }$ has a unique ground state (up to a multiplicative
constant), but $\lambda =0$ is not an eigenvalue. If $d\geq 5$,
then $\lambda =0$ is a simple eigenvalue of $H_{\beta }$ and the
sets of ground states and eigenfunctions coincide.
\end{lemma}
\proof  The ground states belong to $L^{2}(\mathbb{R}^d)$ if and
only if $d\geq 5.$ In order to complete the proof of the lemma, it
remains to show that any eigenfunction of $H_{\beta }$ with zero
eigenvalue satisfies (\ref{grs}). Thus, it is enough to prove that
if $\frac{1}{2}\Delta \psi +\beta v(x)\psi =0$ and $\psi \in
L^{2}( \mathbb{R}^d)$, then $\psi =-R_{0}(0)h$ with $h=\beta
v\psi$. From $\frac{1}{2}\Delta \psi +\beta v(x)\psi -\lambda \psi
=-\lambda \psi$ we obtain $\psi =-R_{0}(\lambda )(h+\lambda \psi
)$. Obviously $R_{0}(\lambda )h\rightarrow R_{0}(0)h$ in
$L^2(\mathbb{R}^d)$ as $\lambda \downarrow 0$ since $h\in
L^2_{\exp}(\mathbb{R}^d)$. Now the lemma  will be proved if we
show that
\[
||\lambda R_{0}(\lambda )\psi
||_{L^{2}(\mathbb{R}^{d})}^{2}=\int_{\mathbb{R}^{d}}
(\frac{2\lambda |\widetilde{\psi}(\sigma)| }{\sigma ^{2}+2\lambda
})^2 d \sigma \rightarrow 0\text{ \ \ \ as \ }\lambda \downarrow
0.
\]
The latter follows from the dominated convergence theorem.\qed

 The following lemma summarizes some facts about the operator
$(I + \beta A(\lambda))^{-1}$ proved above. It also describes
 the structure of the
singularity of the operator $(I + \beta A(\lambda))^{-1}$ for
$\lambda$ and $\beta$ in a neighborhood of $\lambda = 0$, $\beta =
\beta_{cr}$.
\begin{lemma} \label{tet} Let $d \geq 3$ and $\beta \geq 0$.
 The operator $(I + \beta A(\lambda))^{-1}$ (considered in
 $L^2_{\exp}(\mathbb{R}^d)$ and
$C_{\exp}( \mathbb{R}^d)$) is meromorphic in $\lambda \in
\mathbb{C}'$
 and has poles of the first order at
eigenvalues of the operator $H_{\beta}$. For each $\varepsilon
> 0$ and some $\Lambda=\Lambda(\beta)$, the operator is uniformly bounded
in $\lambda \in \mathbb{C}'$,
 $|{\rm arg} \lambda| \leq \pi -
\varepsilon$, $|\lambda| \geq \Lambda$.

If $\beta= \beta_{cr}$, then the operator $(I + \beta
A(\lambda))^{-1}$ is analytic in $\lambda \in \mathbb{C}'$ and
uniformly bounded in $\lambda \in \mathbb{C}'$,
 $|{\rm arg} \lambda| \leq \pi -
\varepsilon$, $|\lambda| \geq \varepsilon$.

 If $\beta <
\beta_{cr}$, then the operator $(I + \beta A(\lambda))^{-1}$ is
analytic in $\lambda \in \mathbb{C}'$ and uniformly bounded in
$\lambda \in \mathbb{C}'$,
 $|{\rm arg} \lambda| \leq \pi -
\varepsilon$.

There are $\lambda_0 > 0$ and $\delta_0 > 0$ such that for
$\lambda \in \mathbb{C}' \cup \{ 0\}$,  $|\lambda| \leq
\lambda_0$, $|\beta - \beta_{cr}| \leq \delta_0$, $\beta \neq
\beta(\lambda)$, we have the representation
\begin{equation} \label{mre}
(I + \beta A(\lambda))^{-1}  =  \frac{\beta(\lambda)
}{\beta(\lambda) - \beta}(B + S_d(\lambda)) + C(\lambda, \beta).
\end{equation}
Here $\beta(\lambda)$ is defined in Remark 2 following
Lemma~\ref{ll44},  $B$ is the one dimensional operator with the
kernel
\begin{equation} \label{bbi}
B(x,y) = \frac{v(x) \psi (x) \psi (y)}{ \int
_{\mathbb{R}^d}v(x)\psi ^2(x)dx},
\end{equation}
where $\psi$ is a ground state defined in the Remark following Lemma \ref{re1}, and
\begin{equation} \label{ssj}
S_3(\lambda)= O(\sqrt{|\lambda|}),~~ S_4(\lambda)= O(|\lambda \ln
(\lambda)|), ~~ S_d(\lambda) = O(|\lambda|),~~d \geq 5,~~{\rm
as}~~ \lambda \rightarrow 0,~~\lambda \in \mathbb{C}',
\end{equation}
$S_d(0) = 0$, $d \geq 3$,
 and $C(\lambda, \beta)$ is bounded
uniformly in $\lambda$ and $\beta$.
\end{lemma}
\proof  The analytic properties of $(I + \beta A(\lambda))^{-1}$
follow from Lemma \ref{re1}. By Lemma~\ref{ker10}, the norm of
$A(\lambda)$ decays at infinity when $\lambda \rightarrow \infty$,
$|{\rm arg} \lambda| \leq \pi - \varepsilon$. Therefore there is
$\Lambda > 0$ such that the operator $(I + \beta A(\lambda))^{-1}$
is bounded for
 $|{\rm arg} \lambda| \leq \pi -
\varepsilon$, $|\lambda| \geq \Lambda$.

If $\beta \leq \beta_{cr}$, then $(I + \beta A(\lambda))^{-1}$
does not have poles in $\lambda \in \mathbb{C}'$, and therefore
$\Lambda$ can be taken to be arbitrarily small.

If $\beta < \beta_{cr}$, then $(I + \beta A(0))$ is invertible by
Lemma~\ref{ll44}. By Lemma~\ref{ker10}, the operators $A(\lambda)$
tend to $A(0)$ when $\lambda \rightarrow 0$, $\lambda \in
\mathbb{C}'$. Therefore $(I + \beta A(\lambda))^{-1}$, $\lambda
\in \mathbb{C}'$, are bounded in a neighborhood of zero. It
remains to justify~(\ref{mre}).

For $d \geq 3$, let $h_\lambda$ be an eigenvector corresponding to
the eigenvalue $1/\beta(\lambda)$ of the operator $-A(\lambda)$,
$\lambda \in [0,\infty) \cup (U \cap \mathbb{C}')$. By
Lemma~\ref{ll44} and the second remark following it, this
eigenvector is defined up to a multiplicative constant.
Let $A^*(\lambda)$ be the operator in $L^2_{\exp}(\mathbb{R}^d)$ or
$C_{\exp}(\mathbb{R}^d)$ with the kernel $A^*(\lambda,x,y) =
A(\lambda,y,x)e^{|y|^2-|x|^2}$.
%$A^*(\lambda)$ be the operator adjoint with respect to the scalar product in
%$L^2(\mathbb{R}^d)$ to the operator $A(\lambda)$, i. e. $A^*(\lambda)$
%is the operator in the space $(L^2_{\exp}(\mathbb{R}^d))^*$ with the kernel $A^*(\lambda,x,y) =
%A(\lambda,y,x)=R_0(\lambda,x-y)v(y)$. The space $(L^2_{\exp}(\mathbb{R}^d))^*$ is adjoint to
%$L^2_{\exp}(\mathbb{R}^d)$, i. e. it is the space $L^2(\mathbb{R}^d)$ with the
%weight $e^{-|x|^2}$.
Similarly to Lemma~\ref{ll44},
it is not difficult to show that $1/\beta(\lambda)$ is an
eigenvalue for the operator $-A^*(\lambda)$ and that its real part
exceeds the absolute values of the other eigenvalues. The
corresponding eigenvector $h^*_\lambda$ is uniquely defined up to
a multiplicative constant. Moreover, we can take $h_\lambda$ and
$h^*_\lambda$ such that
\begin{equation} \label{h*}
v(x)e^{|x|^2}h_\lambda^*(x) =
h_\lambda(x).
\end{equation}
 Note that $h_\lambda$ and $h_\lambda^*$ can be chosen in
such a way that
\begin{equation} \label{hhjj}
 ||h_\lambda - h_0||, ||h^*_\lambda - h^*_0|| \leq
k ||A(\lambda) - A(0)||
\end{equation}
 for some $k > 0$ and all sufficiently
small $|\lambda|$,
where the norms on both sides of (\ref{hhjj})
are either in the space $L^2_{\exp}(\mathbb{R}^d)$ or $C_{\exp}(\mathbb{R}^d)$.

Recall that $A(\lambda) \rightarrow A(0)$ as $\lambda \rightarrow
0$, $\lambda \in \mathbb{C}'$, by Lemma~\ref{ker10}. Using this
and the fact that $1/\beta_{cr}$ is the principal eigenvalue for
$-A(0)$, it is easy to show that there are $\lambda_1 > 0$ and
$\delta_1 > 0$ such that for $\lambda \in \mathbb{C}' \cup \{0\}$,
$|\lambda| \leq \lambda_1$, the eigenvalue $1/\beta(\lambda)$ of
the operator $-A(\lambda)$ is the unique eigenvalue whose distance
from $1/\beta_{cr}$ does not exceed $\delta_1$. Take $0 <
\lambda_0 < \lambda_1$ and $0 < \delta_0 <\delta_1$ such that
 for $\lambda \in \mathbb{C}' \cup \{ 0 \}$, $|\lambda| \leq
\lambda_0$, the distance between  $1/\beta(\lambda)$ and
$1/\beta_{cr}$ does not exceed $\delta_0$.

 Then for $\lambda \in \mathbb{C}' \cup \{0\}$, $|\lambda| \leq
\lambda_0$ and $\beta$ such that $|1/\beta - 1/\beta_{cr}| \leq
\delta_0$, the operator valued function
\[
F(z) = \frac{( A(\lambda) + z I)^{-1}}{ z -({1}/{\beta})}
\]
 is
meromorphic inside the circle $\gamma = \{z: |z - 1/\beta_{cr}| =
\delta_1 \}$. It has two poles: one at $z = 1/\beta$ and the other
at $z = 1/\beta(\lambda)$. The residue at the first pole is equal
to $(A(\lambda) + I/\beta)^{-1}$. In order to find the residue at
the second pole, recall that it is a simple pole for $( A(\lambda)
+ z I)^{-1}$, and therefore
\[
( A(\lambda) + z I)^{-1} = T_{-1}(\lambda)(z -
\frac{1}{\beta(\lambda)})^{-1} + T_0(\lambda) + T_1(\lambda)(z -
\frac{1}{\beta(\lambda)}) + ...
\]
for some operators $T_{-1}, T_0, T_1,...$ and all $z$ in a
neighborhood of $1/\beta(\lambda)$. From here and the fact that
the kernels of $A(\lambda) + I/\beta(\lambda)$ and $A^*(\lambda) +
I/\beta(\lambda)$  are one-dimensional and coincide with ${\rm
span} \{h_\lambda\}$ and ${\rm span} \{h^*_\lambda\}$,
respectively, it easily follows that
\[
T_{-1}(\lambda) f =   \frac{h_\lambda \langle f, h^*_\lambda
\rangle_{L^2_{\exp}(\mathbb{R}^d)}}{\langle h_\lambda ,
h^*_\lambda \rangle_{L^2_{\exp}(\mathbb{R}^d)}},~~f \in
L^2_{\exp}(\mathbb{R}^d)~~({\rm in}~{\rm particular}~{\rm if }~f
\in C_{\exp}(\mathbb{R}^d)).
\]
From (\ref{hhjj}) and  Lemma~\ref{ker10} it follows that
$S_d(\lambda) := T_{-1}(\lambda) - T_{-1}(0)$ satisfies (\ref{ssj}).
The residue of $F(z)$ at $z = 1/\beta(\lambda)$ is
equal to
\[
\frac{\beta(\lambda)\beta }{\beta - \beta(\lambda)}(T_{-1}(0) +
S_d(\lambda)).
\]
Integrating $F(z)$ over the contour $\gamma$, we obtain
\[
( A(\lambda) +  I/\beta)^{-1} + \frac{\beta(\lambda)\beta }{\beta
- \beta(\lambda)}(T_{-1}(0) + S_d(\lambda)) = \frac{1}{2 \pi i}
\int_\gamma \frac{( A(\lambda) + z I)^{-1}}{ z -({1}/{\beta})} d
z.
\]
The right hand side of this formula is uniformly bounded, which
completes the proof of the lemma if we show that $T_{-1}(0)=B$. Thus it remains to prove that
\[
\frac{h_0(x) e^{|y|^2}h^*_0(y) }{\langle h_0 , h^*_0
\rangle_{L^2_{\exp}(\mathbb{R}^d)}}= \frac{v(x) \psi (x) \psi
(y)}{ \int _{\mathbb{R}^d}v(x)\psi ^2(x)dx}.
\]
The latter follows from the relation $h_0=\beta v \psi$ (see Lemma
\ref{re1}) and (\ref{h*}). \qed

Formula (\ref{aa3}) and Lemmas~\ref{conti1} and \ref{tet} imply
the following result.

\begin{lemma} \label{rbeta} Let $d \geq 3$ and $\beta \geq 0$.
 The operator $R_\beta(\lambda)$ (considered as an operator from
$C_{\exp}( \mathbb{R}^d)$ to $C( \mathbb{R}^d)$) is meromorphic in
$\lambda \in \mathbb{C}'$
 and has poles of the first order at
eigenvalues of the operator $H_{\beta}$. For each $\varepsilon
> 0$ and some $\Lambda=\Lambda(\beta)$, the operator is uniformly bounded
in $\lambda \in \mathbb{C}'$,
 $|{\rm arg} \lambda| \leq \pi -
\varepsilon$, $|\lambda| \geq \Lambda$. It is of order
$O(1/|\lambda|)$ as $\lambda \rightarrow \infty$, $|{\rm arg}
\lambda| \leq \pi - \varepsilon$.

If $\beta= \beta_{cr}$, then the operator $R_\beta(\lambda)$ is
analytic in $\lambda \in \mathbb{C}'$ and uniformly bounded in
$\lambda \in \mathbb{C}'$,
 $|{\rm arg} \lambda| \leq \pi -
\varepsilon$, $|\lambda| \geq \varepsilon$.

 If $\beta <
\beta_{cr}$, then the operator $R_\beta(\lambda)$ is analytic in
$\lambda \in \mathbb{C}'$ and uniformly bounded in $\lambda \in
\mathbb{C}'$,
 $|{\rm arg} \lambda| \leq \pi -
\varepsilon$.

There are $\lambda_0 > 0$ and $\delta_0 > 0$ such that for
$\lambda \in \mathbb{C}'$,  $0<|\lambda| \leq \lambda_0$, $|\beta
- \beta_{cr}| \leq \delta_0$, $\beta \neq \beta(\lambda)$, we have
the representation
\begin{equation} \label{mre2}
R_\beta(\lambda)  =  \frac{\beta(\lambda) }{\beta(\lambda) -
\beta}({R}_{0}(0) B + S_d(\lambda)) + C(\lambda, \beta),
\end{equation}
where  $\beta(\lambda)$ is defined in Remark 2 following
Lemma~\ref{ll44} and $B$ is given by (\ref{bbi}), $S_d$, $d \geq
3$, satisfy (\ref{ssj}), and $C(\lambda, \beta)$ is bounded
uniformly in $\lambda$ and $\beta$.

\end{lemma}

\section{The Behavior of the Principal Eigenvalue for $\beta
\downarrow \beta_{cr}$} \label{secff}

 In Lemma~\ref{ll44} we showed that
$\beta_{cr}>0$ for $d \geq 3$. The following theorem implies, in
particular, that $\beta_{cr} = 0$ for $d=1$ or $2$.
\begin{theorem}
For $d =1,2$ (when $\beta_{cr}=0$) the eigenvalue $\lambda
_{0}(\beta )$ has the following behavior as $\beta \downarrow
\beta _{cr}$: \label{tbcr}
\begin{equation}
\lambda_{0}(\beta )\sim \frac{1}{2}c_{1}^{2}\beta
^{2},~~c_{1}=\int_{\mathbb{R}^d }v(x)d x,~~~d=1,  \label{dim1}
\end{equation}
\begin{equation}
\lambda_{0}(\beta )\sim \exp (-\frac{c_{2}}{\beta
}),~~~c_{2}=\frac{\pi }{ c_{1}},~~~d=2.  \label{dim2}
\end{equation}
In dimensions $d \geq 3$ the eigenvalue $\lambda _{0}(\beta )$ has
the following behavior as $\beta \downarrow \beta _{cr}$:
\begin{equation}
\lambda_{0}(\beta )\sim c_{3}(\beta -\beta _{cr})^{2},~~~d=3,
\label{dim3}
\end{equation}
\begin{equation}
\lambda_{0}(\beta )\sim c_{4}(\beta -\beta _{cr})/\ln (1/(\beta
-\beta _{cr})),~~~d=4,  \label{dim4}
\end{equation}
\begin{equation}
\lambda_{0}(\beta )\sim c_d(\beta -\beta _{cr}),~~~d\geq 5,
\label{dim5}
\end{equation}
where $c_{d}\neq 0,$ $d\geq 3,$\ depend on $v$ and will be
indicated in the proof.
\end{theorem}

\proof

%Let $\mu (\lambda )=1/\beta (\lambda )$.
 Since we are interested
in the behavior of $\lambda_0(\beta)$ for $\beta \downarrow
\beta_{cr}$ and $\lambda_0(\beta) \downarrow 0$ when $\beta
\downarrow \beta_{cr}$ by Lemma~\ref{lmonot}, we shall study the
behavior of $\beta(\lambda)$ as~$\lambda \downarrow 0$ (or, more
generally, as $\lambda \rightarrow 0$, $\lambda \in \mathbb{C}'$).
 The arguments below are based on Lemma \ref{ker10}.

First consider the case $d=1$. For  $\lambda \rightarrow 0$,
$\lambda \in \mathbb{C}'$, the eigenvalue problem for $-A(\lambda
)$ can be written in the form
\begin{equation}
(v{P_{1}}+O(\sqrt{{\lambda }}))h_\lambda=\frac{\sqrt{\lambda
}}{\beta (\lambda )}h_\lambda.  \label{yy1}
\end{equation}
Note that the kernel of $vP_{1}$ is positive when $x$ is an
interior point of ${\rm supp} (v) $. Therefore $vP_{1}$ has a
principal eigenvalue. In fact, the operator $vP_{1}$ is
one-dimensional and the eigenvalue is equal to $c_1/\sqrt2$ where $
c_1=\int_{\mathbb{R}^d }v(x)d x$. Since this eigenvalue is
simple and the operator in the left-hand side of (\ref{yy1}) is
analytic in $\sqrt{\lambda } $, both $h_\lambda$ and
$\sqrt{\lambda } /\beta (\lambda )$ are analytic functions of
$\sqrt{\lambda }$ in a neighborhood of the origin and
\[
\lim_{\lambda \rightarrow 0, \lambda \in
\mathbb{C}'}(\sqrt{\lambda }/\beta (\lambda ))=c_{1}/\sqrt2.
\]
Therefore, $\beta _{cr}=0,$ $\beta (\lambda )$\ is analytic in $\sqrt{%
\lambda },$ and $\beta (\lambda )\sim \sqrt{2\lambda }/c_{1}$ as
$\lambda \rightarrow 0$, $\lambda \in \mathbb{C}'$, which proves
(\ref{dim1}).

The same arguments in the case $d=2$ lead to the relation
\[
\lim_{\lambda \rightarrow 0, \lambda \in \mathbb{C}'}
(\frac{-1}{\beta(\lambda) \ln \lambda} ) = c_1/\pi.
\]
This implies that $\beta_{cr} = 0$ and (\ref{dim2}) holds.

In the case $d=3$ the eigenvalue problem for $-A(\lambda )$ takes the form
\begin{equation}
(-A(0)+\sqrt{\lambda }v(x)Q_{3}+O({\lambda }))h_\lambda=\frac{1
}{\beta (\lambda )}h_\lambda. \label{yy3}
\end{equation}
As in the one-dimensional case, $1/\beta (\lambda )$ and
$h_\lambda$ are analytic functions of $\sqrt{\lambda }$. Now
$1/\beta_{cr}$ is equal to the principal eigenvalue of $-A(0)$.
Recall that $h_{0}$ is the principal eigenfunction of $-A(0)$ and
$h_0^*$ is the principal eigenfunction of $-A^*(0)$.  Standard
perturbation arguments imply that
\begin{equation} \label{betas}
\frac{1}{\beta (\lambda )}= \frac{1}{\beta_{cr}}-\gamma
\sqrt{\lambda }+O(\lambda ),~~~ \lambda \rightarrow 0,~ \lambda
\in \mathbb{C}',
\end{equation}
where
\begin{equation}
\gamma =\frac{-\langle vQ_{3}h_{0},h_{0}^{\ast }\rangle _{L^{2}_{\exp}(\mathbb{R}^d)}}{%
\langle h_{0},h_{0}^{\ast }\rangle _{L^{2}_{\exp}(\mathbb{R}^d)}} > 0,
%  =\frac{\langle Q_{3}h_{0},h_{0}\rangle _{L^{2}_{\exp}(\mathbb{R}^d)}}{%
%\langle h_{0},h_{0}^{\ast }\rangle _{L^{2}_{\exp}(\mathbb{R}^d)}}.
\label{gamma}
\end{equation}
which implies (\ref{dim3}) with $c_{3}=1/(\gamma ^{2}\beta
_{cr}^{4})$. Note that $\gamma > 0$ since the kernel of the
operator $vQ_{3}$ is negative and principal eigenfunctions
$h_{0},$ $h_{0}^{\ast }$ can be chosen to be positive inside ${\rm
supp}(v)$.

Formula for $\gamma$ can be simplified. We choose $h_0=\beta v \psi$ (see Lemma \ref{re1}) and $h_0^*$
defined in (\ref{h*}). Then
\begin{equation}
\gamma =\frac{(\int_{\mathbb{R}^3}v(x)\psi(x)dx)^2}{\sqrt{ 2} \pi
\int_{\mathbb{R}^3}v(x)\psi^2(x)dx},~~d=3.
\label{gamma3}
\end{equation}

Let $d=4$. Then instead of  (\ref{yy3}) we get
\begin{equation}
(-A(0)+\lambda \ln (1/\lambda)
vQ_{4}+O({\lambda}))h_\lambda=\frac{1}{\beta(\lambda )}h_\lambda.
\label{yy4}
\end{equation}
From here it follows that
\begin{equation}
\frac{1}{\beta (\lambda )}=\frac{1}{\beta_{cr}}-\gamma \lambda
\ln(1/ \lambda) +O(\lambda ),~~~ \lambda \rightarrow 0,~ \lambda
\in \mathbb{C}', \label{yyy}
\end{equation}
where $1/\beta_{cr}$ is the principal eigenvalue of $-A(0)$ and
$\gamma $ is given by (\ref{gamma}) with $\ Q_{3}$ replaced by $\
Q_{4}.$ Thus
 (\ref{dim4}) holds with $c_{4}=1 / (\gamma\beta _{cr}^{2})$.

For $d \geq 5$ we get
\[
(-A(0)+{\lambda
v}Q_{d}+O({\lambda^{3/2}}))h_\lambda=\frac{1}{\beta(\lambda
)}h_\lambda.
\]
From here it follows that
\[
\frac{1}{\beta (\lambda )}=\frac{1}{\beta_{cr}}-\gamma \lambda
+O(\lambda^{3/2}),~~~ \lambda \rightarrow 0,~ \lambda \in
\mathbb{C}',
\]
where $1/\beta_{cr}$ is the principal eigenvalue of $-A(0)$ and
$\gamma$ is given by (\ref{gamma}) with $\ Q_{3}$ replaced by $\
Q_{d}$. Thus
 (\ref{dim5}) holds with $c_{d}=1 / (\gamma\beta _{cr}^{2})$.
 \qed

\section{Asymptotics of the Partition Function, Solutions,
and Fundamental Solutions} \label{secfi}

We shall need the following notation.  Recall from (\ref{pbeta})
that by $p_\beta(t, y,x)$ we denote the fundamental solution of
the parabolic problem
\[
\frac{\partial p_\beta(t,y,x)}{\partial t} =  \frac{1}{2}\Delta_x
p_\beta(t,y,x) +\beta v(x) p_\beta(t,y,x),
\]
\[
p_\beta(0,y,x) = \delta(x-y).
\]
For a given $f \in L^2( \mathbb{R}^d)$, let
\[
u_\beta(t,x) = \int_{ \mathbb{R}^d}  p_\beta(t,y,x) f(y) d y
\]
 be the
solution of the Cauchy problem with the initial data $f$. The
partition function is defined as the integral of the fundamental
solution
\[
Z_{\beta,t} (x)  = \int_{ \mathbb{R}^d} p_\beta(t,x,y) d y =\int_{
\mathbb{R}^d} p_\beta(t,y,x) d y.
\]
Note that the partition function defined in (\ref{partf}) is
simply $ Z_{\beta,T} = Z_{\beta,T} (0)$. Also note that
$Z_{\beta,t} (x)$ is the solution of the Cauchy problem with
initial data equal to one:
\[
\frac{\partial Z_{\beta,t}(x)}{\partial t} =  \frac{1}{2}\Delta
Z_{\beta,t}(x) +\beta v(x) Z_{\beta,t}(x),~~Z_{\beta,0}(x) \equiv
1.
\]

For $\beta > \beta_{cr}$, let $\psi_\beta$ be the positive
eigenfunction for the operator $H_\beta$ with eigenvalue
$\lambda_0(\beta)$ normalized by the condition
$||\psi_\beta||_{L^2( \mathbb{R})} = 1$. This function is defined
uniquely by Lemma~\ref{re1} and is equal to $-R_0(\lambda)
h_\lambda$, where $\lambda = \lambda_0(\beta)$ and $h_\lambda$ is
the principal eigenfunction for the operator~$-A(\lambda)$. Note
that $\psi_\beta$ decays exponentially at infinity.

 For $a \in \mathbb{R}$,
let $\Gamma(a)$ be the following contour in the complex plane
\[
\Gamma(a) = \{a -s + i s,~s \geq 0\} \cup \{a - s- i s,~s \geq
0\}.
\]
We choose the direction along $\Gamma(a)$ in such a way that the
imaginary coordinate increases.

 The following lemma is an important tool for investigating the
asymptotics of ${Z_{\beta,T}}$.
\begin{lemma} \label{cont1}
Let $a > \lambda_0(\beta)$. Then for $f \in L^2( \mathbb{R}^d)$
(or $f \in C_{\exp}( \mathbb{R}^d)$) and $t
> 0$,
\begin{equation} \label{y5}
u_{\beta}(t,x) = \frac{-1}{2 \pi i} \int_{\Gamma(a)} e^{\lambda t}
(R_\beta(\lambda) f)(x) d \lambda,
\end{equation}
which holds in $L^2( \mathbb{R}^d)$ (or $C( \mathbb{R}^d)$). This
formula remains valid if the initial function $f$ is identically
equal to one and $R_\beta(\lambda) f$ is understood by
substituting $f \equiv 1$ into (\ref{aa1}) with $R_0(\lambda) 1 =
-1/\lambda$. More precisely,
\begin{equation}
 Z_{\beta,t}(x) - 1 = \frac{-1}{2 \pi i} \int_{\Gamma(a)}
\frac{e^{\lambda t}}{\lambda} (R_\beta(\lambda)(\beta v)) (x) d
\lambda \label{yy11}
\end{equation}
%\begin{equation}
% =\frac{1}{2 \pi} \int_{\Gamma(a)} \frac{e^{z t}}{z}
%R_{0}(z )(I+\beta A(z ))^{-1}(\beta v) (x) d z
%\end{equation}
in $L^2( \mathbb{R}^d)$ and $C( \mathbb{R}^d)$.
\end{lemma}
%\noindent {\bf Remark.} The right-hand sides of (\ref{y5}) and
%(\ref{yy11}) are continuous functions of $x$, as will be seen in
%the proof of Theorem~\ref{thz}.
\proof First, let $f \in L^2( \mathbb{R}^d)$. We solve the Cauchy
problem for $u_\beta$ using the Laplace transform with respect
to~$t$. This leads to (\ref{y5}) with $\Gamma(a)$ replaced by the
line $\{\lambda : {\rm Re}\lambda = a \}$. The integral over this
line is equal to the integral over $\Gamma(a)$ since the resolvent
is analytic between these contours and its norm decays as
$|\lambda|^{-1}$ when $|\lambda| \rightarrow \infty$.

Now let $f \equiv 1$. Then $w(t,x) = Z_{\beta, t}(x) - 1$ is the
solution of the problem
\[
\frac{\partial w(t,x)}{\partial t} = \frac{1}{2}\Delta w(t,x) +
\beta v(x) w(t,x) + \beta v(x),~~~w(0,x) \equiv 0.
\]
By the Duhamel formula and (\ref{y5}),
\[
 w(t,x) = \frac{-1}{2 \pi i} \int_0^t \int_{\Gamma(a)} e^{\lambda( t-s)}
(R_\beta(\lambda) \beta v)(x) d \lambda d s =
\]
\[
\frac{-1}{2 \pi i}\int_{\Gamma(a)} \frac{e^{\lambda t}-1}{\lambda}
(R_\beta(\lambda) \beta v)(x) d \lambda = \frac{-1}{2 \pi
i}\int_{\Gamma(a)} \frac{e^{\lambda t}}{\lambda} (R_\beta(\lambda)
\beta v)(x) d \lambda,
\]
since in the domain $\Gamma^+(a)$ to the right of the contour
$\Gamma(a)$, the operator $ R_\beta(\lambda): L^2( \mathbb{R}^d)
\rightarrow L^2( \mathbb{R}^d)$ is analytic and decays as
$|\lambda|^{-1}$ at infinity. This justifies (\ref{yy11}) in $L^2(
\mathbb{R}^d)$ sense. It remains to show that the right-hand side
of (\ref{y5}) is continuous for $f \in C_{\exp}( \mathbb{R}^d)$
and the right-hand side of (\ref{yy11}) is continuous. Since
$\beta v \in C_0^\infty$, the integrands are continuous  in
$(t,x)$ for each $\lambda \in \Gamma_a$. It remains to note that
the integrals converge uniformly when $x \in \mathbb{R}^n$, $t
\geq t_0>0$. This is due to the fact that $ ||R_\beta(\lambda)
f||_{C( \mathbb{R}^d)}, ||R_\beta(\lambda) \beta v||_{C(
\mathbb{R}^d)} \leq C_d(a)$, as follows from Lemma~\ref{rbeta}.
 \qed

In order to state the next theorem we shall need the following
notation. As in part (3) of Lemma~\ref{re1}, it is not difficult
to show that for $d \geq 3$, $0 \leq \beta < \beta_{cr}$ and $f
\in C_0^\infty( \mathbb{R}^d)$ there is a unique solution of the
problem
\begin{equation} \label{fibeta}
H_{\beta }(\varphi)=\frac{1}{2}\Delta \varphi +\beta v(x)\varphi
=f,~~~\varphi = O(|x|^{2-d}),~~\frac{\partial{\varphi}
}{\partial{r}}(x) = O(|x|^{1-d})~~as~~r =|x| \rightarrow \infty.
\end{equation}
This solution is given by $\varphi = R_0(0)(I + \beta A(0))^{-1}
f$. For $f  = -\beta v$, we denote this solution by
$\varphi_\beta$.
\begin{theorem}
\label{thz} (1) For $\beta > \beta_{cr}$ there is  $\varepsilon >
0$ such that
 we have the following asymptotics for
the partition function:
\[
 Z_{\beta,t}(x)-1 =  \exp(\lambda_0(\beta) t)
( ||\psi_\beta||_{L^1( \mathbb{R}^d)} \psi_\beta(x) +
O(\exp(-\varepsilon t)))~~as~~t \rightarrow
 \infty,
\]
which holds in $L^2( \mathbb{R}^d)$ and  in $C(\mathbb{R}^d)$,
where $\psi_\beta$ is the positive eigenfunction for the operator
$H_\beta$ with eigenvalue $\lambda_0(\beta)$ normalized by the
condition $||\psi_\beta||_{L^2( \mathbb{R})} = 1$.

(2) For $\beta = \beta_{cr}$ we have the following asymptotics for
the partition function:
\[
 Z_{\beta,t}(x) =  k_3 t^{1/2} \psi(x) + O(1)~~as~~t \rightarrow
 \infty,~~d =3,
\]
\[
 Z_{\beta,t}(x) = k_4 \frac{t}{\ln t}\psi(x) + O(\frac{t}{\ln^2 t})  ~~as~~t \rightarrow
 \infty,~~d =4,
\]
\[
 Z_{\beta,t}(x) =  k_d t \psi(x) + O (\sqrt{t})~~as~~t \rightarrow
 \infty,~~d \geq 5,
\]
which holds in $C(\mathbb{R}^d)$. Here $k_d$, $d \geq 3$, are
positive constants and $\psi$ is the positive ground state for
$H_{\beta_{cr}}$ normalized by the condition ${||\beta_{cr} v
\psi||_{L^2_{\exp}(\mathbb{R}^d)} = 1}$.

(3) If $0 \leq \beta < \beta_{cr}$, then
\[
\lim_{t \rightarrow \infty} Z_{\beta,t}(x)=  1+ \varphi_\beta(x)
\]
in $C(\mathbb{R}^d)$.
\end{theorem}
\proof (1) Note that the resolvent $R_\beta(\lambda)$ has only one
pole between the contours $\Gamma(a)$ and $\Gamma(\lambda_0(\beta)
- \varepsilon)$ if $\varepsilon$ is less than the distance from
$\lambda_0$ to the rest of the spectrum. This pole is at the point
$\lambda_0(\beta)$ and the residue is the integral operator with
the kernel $-\psi_\beta(x) \psi_\beta(y)$. Therefore from
(\ref{yy11})  it follows that
\begin{equation} \label{ab1}
Z_{\beta,t}(x) - 1 = \frac{e^{\lambda_0(\beta) t}}{
\lambda_0(\beta)} \psi_\beta(x) \int_{ \mathbb{R}^d} \beta v(y)
\psi_\beta (y) d y- \frac{1}{2 \pi i}
\int_{\Gamma(\lambda_0(\beta) - \varepsilon)} \frac{e^{\lambda
t}}{\lambda} (R_\beta(\lambda) \beta v)(x) d \lambda.
\end{equation}
Since $(\frac{1}{2}\Delta + \beta v - \lambda_0(\beta)) \psi_\beta
= 0$, we have $\beta v \psi_\beta = (\lambda_0(\beta) -
\frac{1}{2}\Delta) \psi_\beta$, and the integral in the first term
of the right-hand side of (\ref{ab1}) is equal to
$\lambda_0(\beta) ||\psi_\beta||_{L^1( \mathbb{R}^d)}$. Thus the
first term on the right-hand side coincides with the main term of
the asymptotics stated in the theorem.

It remains to show that the second term on the right-hand side of
(\ref{ab1}) is exponentially smaller than the first term. This is
due to the fact that the norm of the operator $R_\beta(\lambda)$
is of order $1/|\lambda|$ at infinity for $\lambda \in
\Gamma(\lambda_0(\beta) - \varepsilon)$.

(2) Let $d= 3$. First, let us analyze (\ref{mre2}) when $\beta =
\beta_{cr}$ and $\lambda \rightarrow 0$, $\lambda \in
\mathbb{C}'$. By (\ref{betas}), the factor ${\beta(\lambda)
}/{(\beta(\lambda) - \beta)}$ in the right hand side of
(\ref{mre2}) is equal to $(\beta_{cr}\gamma \sqrt{\lambda})^{-1} +
O(1)$ as $\lambda \rightarrow 0$, $\lambda \in \mathbb{C}'$, where
$\gamma > 0$ is given by (\ref{gamma}).

We choose the same ground state $\psi$  specified in the statement
of Theorem \ref{thz}. Then from (\ref{bbi}) and Lemma~\ref{re1} it
follows that
\begin{equation} \label{r0}
{R}_{0}(0) B (\beta_{cr} v) = \frac{\int_{R^d}v(x) \psi (x) dx}{ \int _{R^d}v(x)\psi ^2(x)dx}
{R}_{0}(0) (\beta_{cr}v\psi) = -\frac{\int_{R^d}v(x) \psi (x) dx}{ \int _{R^d}v(x)\psi ^2(x)dx} \psi.
\end{equation}
%We choose $h_0$ in (\ref{bbi}) (and in the formula above) to be nonnegative and
%normalized
%by the condition $||h_0||_{L^{2}_{\exp}(\mathbb{R}^d)} = 1$. Then the ground state
%defined in the statement of the theorem is equal to $-R_{0}(0) h_0$ (see Lemma~\ref{re1}) and
%this justifies the last equality above.
Now, by Lemma~\ref{rbeta} and (\ref{dim3}), (\ref{betas}),
\begin{equation} \label{cpd}
R_{\beta_{cr}}(\lambda) (\beta_{cr} v ) =
\frac{-\int_{R^d}v(x) \psi (x) dx}{ \gamma \beta_{cr}
\sqrt{\lambda}\int _{R^d}v(x)\psi ^2(x)dx} \psi + D(\lambda)=
\frac{-k'_3 \psi}{\sqrt{\lambda}} +D(\lambda), ~ ~ ~k'_3>0,
\end{equation}
where the remainder $D(\lambda)$ is of order $O(1)$ when $\lambda
\rightarrow 0$, $\lambda \in \mathbb{C}'$. Note that $D(\lambda)$
is bounded on $\Gamma^+(0)$ since the left hand side and the first
term on the right hand side of (\ref{cpd}) are bounded on
$\Gamma^+(0)$ outside a neighborhood of zero.

Next, we apply (\ref{yy11}) with $a$ replaced by $1/t$ and use the
expression (\ref{cpd}) to obtain
\begin{equation} \label{cpd1}
 Z_{\beta,t}(x) - 1 = \frac{1}{2 \pi i} \int_{\Gamma(1/t)}
\frac{e^{\lambda t}}{\lambda} ( \frac{k'_3 \psi}{\sqrt{\lambda}}
+ D(\lambda)) d \lambda.
\end{equation}
Let us change the variables in the integral $\lambda t = z$. Thus
\[
Z_{\beta,t}(x) - 1 = \frac{1}{2 \pi i} \int_{\Gamma(1)} \frac{
e^{z}}{z }
( \frac{ \sqrt{t} k'_3 \psi }{\sqrt{z}}  +
D( \frac{z}{t})) d z.
\]
The contribution to the integral from the term containing
$D({z}/{t})$ is bounded, while the contribution from the first
term is equal to $k_3 t^{1/2} \psi(x)$, as claimed in the lemma.
One needs only to note that $k_3 >0$ since
\[
\frac{1}{2 \pi i} \int_{\Gamma(1)} z^{-3/2} e^{z} d z=
\frac{1}{\pi i} \int_{\Gamma(1)} z^{-1/2} e^{z} d z=\frac
{2}{\pi} \int_0^{\infty} \sigma ^{-1/2}e^{-\sigma } d \sigma =\frac {2}{\sqrt{\pi}}>0.
\]
If $d =4$, then (\ref{dim4}), (\ref{yyy}) imply that $\beta
(\lambda)-\beta_{cr}\sim  \beta^2_{cr}\gamma\lambda
\ln(1/\lambda)$ as $\lambda \rightarrow 0, \lambda \in C'.$ This
leads to the following analog of (\ref{cpd1})
%\[
% Z_{\beta,t}(x) - 1 = \frac{-1}{2 \pi i} \int_{\Gamma(1/t)}
%\frac{e^{\lambda t}}{\lambda} ( \frac{ \beta^2_{cr} \langle h^*_0,
%v \rangle_{L^{2}_{\exp}(\mathbb{R}^d)} }{\gamma {\lambda} \ln
%\lambda \langle h_0 , h^*_0 \rangle_{L^{2}_{\exp}(\mathbb{R}^d)}}
%\psi(x)  + D(\lambda)) d \lambda,
%\]
\[
 Z_{\beta,t}(x) - 1 = \frac{1}{2 \pi i} \int_{\Gamma(1/t)}
\frac{e^{\lambda t}}{\lambda} ( \frac{ k'_4 \psi(x)}{\lambda
\ln(1/ \lambda) }
  + D(\lambda)) d \lambda,~ ~ ~k'_4 > 0,
\]
where  $D(\lambda)$ is of order $O(1/|\lambda \ln^2 \lambda|)$
when $\lambda \rightarrow 0$, $\lambda \in \mathbb{C}'$ and is
bounded at infinity. After the change of variables $\lambda t =
z$, we obtain
%\[
%Z_{\beta,t}(x) - 1 = \frac{1}{2 \pi i} \int_{\Gamma(1)} \frac{
%e^{z}}{{z}  } ( \frac{ {t}\beta^2_{cr} \langle h^*_0, v
%\rangle_{L^{2}_{\exp}(\mathbb{R}^d)} }{(\ln{z} - \ln {t})\gamma
%\langle h_0 , h^*_0 \rangle_{L^{2}_{\exp}(\mathbb{R}^d)}} \psi(x)
%+ D(\frac{z}{t})d z,
%\]
\[
Z_{\beta,t}(x) - 1 = \frac{1}{2 \pi i} \int_{\Gamma(1)} \frac{
e^{z}}{{z}  } ( \frac {t k'_4 \psi (x)}{z(\ln{t} - \ln {z})} +
D(\frac{z}{t}))d z,
\]
which easily leads to the second part of the lemma in the case $d
=4$. The treatment of the case $d \geq 5$ is similar.

(3) We apply (\ref{yy11}) with $a$ replaced by $1/t$ to obtain
\begin{equation} \label{koo}
 Z_{\beta,t}(x) - 1 = \frac{-1}{2 \pi i} \int_{\Gamma(1/t)}
\frac{e^{\lambda t}}{\lambda} (R_\beta(\lambda)(\beta v)) (x) d
\lambda = \frac{-1}{2 \pi i} \int_{\Gamma(1)} \frac{e^{z}}{z}
(R_\beta(\frac{z}{t})(\beta v)) (x) d z.
\end{equation}
Note that by Lemma~\ref{conti1} and since $1/\beta$ is not an
eigenvalue of $A(0)$ we have
\[
\lim_{\lambda \rightarrow 0, \lambda \in \mathbb{C}'}
R_\beta(\lambda)(\beta v) = \lim_{\lambda \rightarrow 0, \lambda
\in \mathbb{C}'} R_{0}(\lambda )(I+\beta A(\lambda ))^{-1} (\beta
v) = R_0(0)(I + \beta A(0))^{-1}  (\beta v) = -\varphi_\beta.
\]
Since the difference between  $R_\beta(z/t)(\beta v)$ and
$-\varphi_\beta$ is bounded on $\Gamma(1)$, one can pass to the
limit $t \rightarrow \infty$ under the integral sign in
(\ref{koo}), which leads to
\[
\lim_{t \rightarrow \infty} Z_{\beta,t}(x) = 1 +
\frac{\varphi_\beta(x)}{2 \pi i} \int_{\Gamma(1)} \frac{e^{z}}{z}
d z = 1 + \varphi_\beta(x).
\]
\qed

The third part of Theorem~\ref{thz} establishes the existence of $
\lim_{t \rightarrow \infty} Z_{\beta,t}(x)$ for $\beta <
\beta_{cr}$. Next we examine the behavior of this quantity as
$\beta \uparrow \beta_{cr}$.
\begin{lemma}
There are positive constant $b_d$, $d \geq 3$, such that
\[
\lim_{t \rightarrow \infty} Z_{\beta,t}(x)-1 =
\frac{b_d}{\beta_{cr} - \beta} \psi(x) + O (1)~~as~~\beta \uparrow
\beta_{cr}
\]
is valid in $C( \mathbb{R}^d)$, where $\psi$ is the positive
ground state for $H_{\beta_{cr}}$ normalized by the condition
${||\beta_{cr} v \psi||_{L^{2}_{\exp}(\mathbb{R}^d)} = 1}$.
\end{lemma}
\proof By the third part of Theorem~\ref{thz}, we only need to
find the asymptotics as $\beta \uparrow \beta_{cr}$ of
$\varphi_\beta =  -R_0(0)(I + \beta A(0))^{-1}(\beta v)$. From
(\ref{mre}) with $\lambda = 0$ and $\beta(0) = \beta_{cr}$ and (\ref{r0})
it follows that
\[
\varphi_\beta =  -R_0(0)(I + \beta A(0))^{-1}(\beta v) =
\frac{-\beta_{cr}}{\beta_{cr} - \beta} R_0(0) B (\beta_{cr} v) +
O(1) = \frac{b_d}{ \beta_{cr} - \beta} \psi + O (1)
\]
for some positive constant $b_d$. \qed

\section{Behavior of the Polymer for  $\protect\beta > \protect\beta_{cr}$}
\label{bp}

In this section we shall assume that $\protect\beta >
\protect\beta_{cr}$ is fixed. A result similar to the first part
of Theorem~\ref{thz} is valid for the solution of the Cauchy
problem and for the fundamental solution.
\begin{theorem} \label{lethz} Let $f \in L^2( \mathbb{R}^d)$
(or $f \in C_{\exp}( \mathbb{R}^d)$). For $\beta > \beta_{cr}$
there is  $\varepsilon > 0$ such that we have the following
asymptotics for the solution $u_{\beta}$ of the Cauchy problem
with the initial data $f$:
\begin{equation} \label{eqeq}
 u_{\beta}(t) =  \exp(\lambda_0(\beta) t)
( \langle \psi_\beta,f\rangle_{L^2( \mathbb{R}^d)} \psi_\beta +
q_f( t)),
\end{equation}
which holds in $L^2( \mathbb{R}^d)$ (or in $C(\mathbb{R}^d)$),
where $||q_f( t)|| \leq c ||f|| \exp(-\varepsilon t)$ for some $c$
and all sufficiently large $t$.

We have the following asymptotics for the fundamental solution of
the parabolic equation:
\begin{equation} \label{fsfs}
p_{\beta}(t,y,x) =  \exp(\lambda_0(\beta) t) ( \psi_\beta(y)
\psi_\beta(x) + q(t,y,x)),
\end{equation}
 where $\lim_{t \rightarrow
\infty} ||q(t,y,x)|| = 0$, uniformly in $y$, and (\ref{fsfs})
holds in $L^2( \mathbb{R}^d)$ and in $C(\mathbb{R}^d)$ for each
$y$ fixed.
\end{theorem}
\proof The proof of (\ref{eqeq}) is the same as the proof of the
first part of Theorem~\ref{thz}, and therefore we omit it.

Let $f_\beta^{\delta,y}(x)= p_\beta(\delta,y,x)$ be the
fundamental solution of the parabolic problem at time~$\delta$.
Note that $f_\beta^{\delta,y} \in L^2( \mathbb{R}^d)$ for all
$\delta
> 0$ and all $y$, and $f_\beta^{\delta,y} \in C_{\exp}( \mathbb{R}^d)$
for all sufficiently small
$\delta > 0$ and all $y$.  Denote the solution of the parabolic
equation with the initial data $f_\beta^{\delta,y}$ by
$u^{\delta,y}_\beta(t,x)$. Then
\[
p_\beta(t,y,x) = u^{\delta,y}_\beta(t-\delta,x) =
\exp(\lambda_0(\beta)( t - \delta)) (
\langle\psi_\beta,f_\beta^{\delta,y}\rangle_{L^2( \mathbb{R}^d)}
\psi_\beta(x) + q^\delta(t,y,x)),
\]
where $||q^\delta(t,y,x)||  \leq c ||f_\beta^{\delta,y}||
\exp(-\varepsilon (t-\delta))$ for some $c$ and all sufficiently
large $t$.

 Note that $\langle\psi_\beta,f_\beta^{\delta,y}\rangle_{L^2(
\mathbb{R}^d)}$ can be made arbitrarily close to $\psi_{\beta}(y)$
uniformly in $y$, by choosing a sufficiently small $\delta$,  and
$||f_\beta^{\delta,y}||$ is uniformly bounded in $y$ for any fixed
$\delta$ . This justifies (\ref{fsfs}). \qed

Next, let us study the distribution of the end of the polymer with
respect to the measure $\mathrm{P}_{\beta,T}$ as $T \rightarrow
\infty$.
% Below we shall repeatedly use the fact that $p_\beta(t,x,y) = p_\beta(t, y,x)$.
\begin{theorem}
The distribution of $x(T)$ with respect to the measure
$\mathrm{P}_{\beta,T}$ converges, weakly, as $T \rightarrow
\infty$, to the distribution with the density $
\psi_\beta/||\psi_\beta||_{L^1( \mathbb{R}^d)}$.
\end{theorem}
\proof The density of $x(T)$ with respect to the Lebesgue measure
is equal to
\begin{equation} \label{hh01}
\frac{p_\beta(T, 0,x)}{Z_{\beta, T}(0)} =
\frac{\exp(\lambda_0(\beta) T) ( \psi_\beta(0) \psi_\beta(x)  +
q(T,0,x))} {\exp(\lambda_0(\beta) T) ( ||\psi_\beta||_{L^1(
\mathbb{R}^d)} \psi_\beta(0) + o(1))},
\end{equation}
 where $q$ is the same as
in (\ref{fsfs}). When $T \rightarrow \infty$, the right hand side
of (\ref{hh01}) converges to $\psi_\beta(x)/||\psi_\beta||_{L^1(
\mathbb{R}^d)}$ uniformly in $x$ by Theorem~\ref{lethz}. This
justifies the weak convergence. \qed

Now let us examine the behavior of the polymer in a region
separated both from zero and $T$. Let $S(T)$ be such that
\begin{equation} \label{sttt}
\lim_{T \rightarrow \infty} S(T) = \lim_{T \rightarrow \infty} (T
- S(T)) = +\infty.
\end{equation}
Let $s > 0$ be fixed. Consider the process $y^T(t) = x(S(T)+t)$,
$0 \leq t \leq s$.
\begin{theorem} \label{pro1}
The distribution of the process $y^T(t)$ with
 respect to either of the measures  $\mathrm{P}_{\beta,T}$ or $\mathrm{P}_{\beta,T}(\cdot|x(T)=0)$
 converges as $T
\rightarrow \infty$, weakly in the space $C([0,s], \mathbb{R}^d)$,
to the distribution of a stationary Markov process with invariant
density $\psi_\beta^2$ and the generator
\[
L_\beta g =  \frac{1}{2}\Delta g + \frac{(\nabla \psi_\beta,
\nabla g)}{\psi_\beta}.
\]
\end{theorem}

\noindent {\bf Remark.} Let
\begin{equation} \label{rprp}
 r_\beta(t,y,x) = \frac{p_\beta(t, y,x)
\psi_\beta(x)}{\psi_\beta(y)} \exp(-\lambda_0(\beta) t).
\end{equation}
Note that $r_\beta(t,y,x)$ is the fundamental solution  for the
operator ${\partial}/{\partial t} - L^*_\beta$, where $L^*_\beta$
is the formal adjoint to $L_\beta$. Thus $ r_\beta$ is the
transition density for the Markov process with the generator
$L_\beta$. Also note that $L^*_\beta \psi_\beta^2 = 0$, and thus
$\psi_\beta^2$ is the invariant density for the Markov process.
% Note that the operators $L_\beta$ and
%$H_\beta$ are related by the formula
%\[
%(L_\beta + \lambda_0(\beta))  g = \psi_\beta^{-1} H_\beta (g
%\psi_\beta).
%\]
%The fundamental solution $r_\beta(t,y,x)$ for the operator
%${\partial}/{\partial t} - L_\beta$ is related to the fundamental
%solution $p_\beta(t,y,x)$ for the operator ${\partial}/{\partial
%t} - H_\beta$ via
%\begin{equation} \label{rprp}
% r_\beta(t,y,x) = \frac{p_\beta(t, y,x)
%\psi_\beta(x)}{\psi_\beta(y)} \exp(-\lambda_0(\beta) t).
%\end{equation}
\\
\\
{\it Proof of Theorem~\ref{pro1}.}  We shall only consider the
measure $\mathrm{P}_{\beta,T}$ since the arguments for the measure
$\mathrm{P}_{\beta,T}(\cdot|x(T)=0)$ are completely analogous.
First, let us prove the convergence of the finite-dimensional
distributions. For $y \in \mathbb{R}^d$ and a Borel set $A \in
\mathcal{B}( \mathbb{R}^d)$, let
\[
R(t,y, A) = \int_{A} r_\beta(t,y,x) d x,
\]
with $r_\beta$ given by (\ref{rprp}). Note that $R$ is a Markov
transition function since
\[
\int_{\mathbb{R}^d} r_\beta(t,y,x) d x \equiv 1.
\]
The generator of the corresponding Markov process is $L_\beta$ and
the invariant density is $\psi_\beta^2$.
  Let $0 \leq t_1 < ... < t_n \leq s$. The density of the random
vector $(y^T(t_1),...,y^T(t_n))$ with respect to the Lebesgue
measure on $ \mathbb{R}^{dn}$ is equal to
\[
\rho^T(x_1,...,x_n) =
\]
\[
{ p_\beta(S(T) + t_1,  0, x_1) p_\beta(t_2- t_1, x_1, x_2)...
p_\beta(t_n - t_{n-1},x_{n-1}, x_{n}) Z_{\beta, T -
t_n}(x_n)}{(Z_{\beta, T}(0))^{-1}}.
\]
We replace here all factors $p_{\beta}$, except the first one, by
$r_{\beta}$ using (\ref{rprp}). We replace the first factor and
the factors $Z$ by their asymptotic expansions given in
 Theorems~\ref{lethz} and~\ref{thz}, respectively. This leads to
%\[
%e^{\lambda_0(\beta)(t_n - t_1)} p_\beta(S(T) + t_1,  0, x_1)
%\psi_\beta(x_1) r_\beta(t_2- t_1, x_1, x_2)...
%\]
%\[
%r_\beta(t_n - t_{n-1},x_{n-1}, x_{n})(\psi_\beta(x_n))^{-1}
%Z_{\beta, T - t_n}(x_n) {(Z_{\beta, T}(0))^{-1} } =
%\]
%\[
%(\psi_\beta(0)\psi_\beta(x_1) + q(S(T) + t_1,0,x_1))
%\psi_\beta(x_1) r_\beta(t_2- t_1, x_1, x_2)...
%\]
%\[
%r_\beta(t_n - t_{n-1},x_{n-1}, x_{n})(\psi_\beta(x_n))^{-1} (
%||\psi_\beta||_{L^1( \mathbb{R}^d)} \psi_\beta(x_n) + o(1)) {(
%||\psi_\beta||_{L^1( \mathbb{R}^d)} \psi_\beta(0) + o(1))^{-1} },
%\]
%due to Theorems~\ref{thz} and \ref{lethz}. When $T \rightarrow
%\infty$, the right hand side converges, uniformly in
%$(x_1,...,x_n)$, to
\[
\rho^T(x_1,...,x_n) = \psi^2_\beta(x_1) r_\beta(t_2- t_1, x_1,
x_2)... r_\beta(t_n - t_{n-1},x_{n-1}, x_{n})+o(1),~~~T\rightarrow
\infty,
\]
where the remainder tends to zero uniformly in $(x_1,...,x_n)$. By
the remark made after the statement of the theorem, this justifies
the convergence of the finite dimensional distributions of $y^T$
to those of the Markov process. It remains to justify the
tightness of the family of measures induced by the processes
$y^T$.

%
%Since $\psi_\beta$ is an eigenfunction of the operator $H_{\beta}$
%with the eigenvalue $\lambda_0(\beta)$, from (\ref{rprp})
% it follows that
%\[
%\int r_{\beta}(t,y,x)\psi^2_{\beta}(y)dy=\psi^2_{\beta}(x),~~\int r_{\beta}(t,y,x)dx=1,
%\]
% which together with the previous formula justifies the convergence of the finite dimensional
%distributions of $y^T$ to those of the Markov process. It remains
%to justify the tightness of the family of measures induced by the
%processes $y^T$.

 From the convergence of the one-dimensional
distributions it follows that for any $\eta
> 0$ there is $a > 0$ such that
\begin{equation} \label{tight1}
\mathrm{P}_{\beta,T} (|y^T(0)| > a) \leq \eta.
\end{equation}
for all sufficiently large $T$.  For a
 continuous function $x: [0,T]
\rightarrow {\mathbb{R}}^d$, $x(0) = 0$, let
\[
m^T(x, \delta) = \sup_{|t_1 - t_2| \leq \delta,~ S(T) \leq t_1,
t_2 \leq S(T) + s} |x(t_1) - x(t_2)|.
\]
Let us prove that for each $\varepsilon, \eta > 0$, there is
$\delta
> 0$ such that
\begin{equation} \label{tight2}
\mathrm{P}_{\beta,T}(m^T(x, \delta) > \varepsilon) \leq \eta
\end{equation}
for all sufficiently large $T$. Observe that
\[
\mathrm{P}_{\beta,T}(m^T(x, \delta) > \varepsilon)  =
\]
\[
(Z_{\beta, T}(0))^{-1} \mathrm{E}_{0,T}(\exp(\int_0^{S(T) + s}
\beta v (x(t)) d t) \chi_{\{m^T(x, \delta) > \varepsilon\}}
Z_{\beta, T-S(T)-s}(x(S(T)+s))) \leq
\]
\[
(Z_{\beta, T}(0))^{-1} \sup_{x \in \mathbb{R}^d} Z_{\beta,
T-S(T)-s}(x) \mathrm{E}_{0,T}(\exp(\int_0^{S(T) + s} \beta v
(x(t)) d t) \chi_{\{m^T(x, \delta) > \varepsilon\}} ) \leq
\]
\[
\exp(s \beta \sup_{x \in \mathbb{R}^d} v(x))(Z_{\beta, T}(0))^{-1}
\sup_{x \in \mathbb{R}^d} Z_{\beta, T-S(T)-s}(x)
\mathrm{E}_{0,T}(\exp(\int_0^{S(T)} \beta v (x(t)) d t)
\chi_{\{m^T(x, \delta) > \varepsilon\}} )
\]
\[
\leq \exp(s \beta \sup_{x \in \mathbb{R}^d} v(x))(Z_{\beta,
T}(0))^{-1} \sup_{x \in \mathbb{R}^d} Z_{\beta, T-S(T)-s}(x)
\sup_{x \in \mathbb{R}^d} p_\beta(S(T), 0, x) C(\delta,
\varepsilon),
\]
where $C(\delta, \varepsilon)$ is the probability that for a
$d$-dimensional Brownian motion $W_t$, $0 \leq t \leq s$, we have
\[
\sup_{|t_1 - t_2| \leq \delta,~ 0 \leq t_1, t_2 \leq s} |W(t_1) -
W(t_2)| > \varepsilon.
\]
Note that
\[
\exp(s \beta \sup_{x \in \mathbb{R}^d} v(x))(Z_{\beta, T}(0))^{-1}
\sup_{x \in \mathbb{R}^d} Z_{\beta, T-S(T)-s}(x) \sup_{x \in
\mathbb{R}^d} p_\beta(S(T), 0, x)
\]
is bounded, as follows from Theorems~\ref{thz} and \ref{lethz},
while $C(\delta, \varepsilon)$ can be made arbitrarily small by
selecting a sufficiently small $\delta$. This justifies
(\ref{tight2}). Since the inequalities (\ref{tight1}) and
(\ref{tight2}) hold for all sufficiently large $T$, by choosing
different $a$ and $\delta$, we can make sure that they hold for
all $T$. Thus the family of measures induced by the processes
$y^T$ is tight.
 \qed
\\

\noindent
 {\bf Remark.} If instead of (\ref{sttt}) we assume
that $S(T) = 0$, the result of Theorem~\ref{pro1} will hold with
the only difference that the initial distribution for the limiting
Markov process will now be concentrated at zero, instead of being
the invariant distribution.

\section{Behavior of the Polymer for $\beta < \beta_{cr}$}
\label{bb1x}

%In this section we shall restrict ourselves to the case $d = 3$.
%Results similar to Lemma~\ref{uu1} and Theorem~\ref{nn1}, below,
%are valid in higher dimensions, but the proofs are more
%complicated.

 First,
we shall study the asymptotic behavior of the solution
$u_\beta(t,x)$ of the Cauchy problem and of the fundamental
solution $p_\beta(t,y,x)$ when $t \rightarrow \infty$, $|y| \leq
\varepsilon^{-1}$, $\varepsilon \sqrt{t} \leq |x| \leq
\varepsilon^{-1} \sqrt{t}$, and $\varepsilon > 0$ is small but
fixed. Recall that $\varphi_\beta$ was defined before
Theorem~\ref{thz}.
\begin{lemma} \label{uu1}
Let $d \geq 3$,  $0 \leq \beta < \beta_{cr}$, $\varepsilon > 0$
and $f \in C_{\exp}( \mathbb{R}^d)$, $f \geq 0$. We have the
following asymptotics for the solution $u_{\beta}$ of the Cauchy
problem with the initial data $f$:
\begin{equation} \label{euu0}
u_{\beta}(t,x) =  (2 \pi t)^{-d/2} \exp(-|x|^2/2t) ( \langle 1 +
\varphi_\beta , f \rangle_{L^2( \mathbb{R}^3)}  + q_f(t,x) ),
\end{equation}
where for some constant $C_\beta(\varepsilon)$ we have
\[
 \sup_{\varepsilon \sqrt{t} \leq |x|
\leq \varepsilon^{-1} \sqrt{t}}| q_f(t,x)| \leq
C_\beta(\varepsilon) t^{-1/2} ||f||_{ C_{\exp}( \mathbb{R}^3)},~~t
\geq 1.
\]

We have the following asymptotics for the fundamental solution of
the parabolic equation:
\begin{equation} \label{euu1}
p_\beta(t,y,x) =  (2 \pi t)^{-d/2} \exp(-|x|^2/2t)( 1 +
\varphi_\beta(y)+ q(t,y,x)),
\end{equation}
where
\[
 \lim_{t \rightarrow \infty} \sup_{|y| \leq \varepsilon^{-1},~
\varepsilon \sqrt{t} \leq |x| \leq \varepsilon^{-1} \sqrt{t}}
|q(t,y,x| = 0.
\]
\end{lemma}
\proof Note that (\ref{euu1}) follows from (\ref{euu0}) since the
fundamental solution at time $t$ is equal to the solution with the
initial data $p_\beta(t, y,\delta)$ evaluated at time $t - \delta$
(the same argument was used in the proof of Theorem~\ref{lethz}).
Therefore it is sufficient to prove (\ref{euu0}).

For the sake of transparency of exposition, we shall consider only
the case $d=3$. From Lemma~\ref{rbeta} it follows that we can put
$a=0$ in (\ref{y5}) when $\beta <\beta _{cr}$. Thus
using~(\ref{aa3}) and the explicit formula for $R_0(\lambda)$, we
obtain
\begin{equation}
u_{\beta }(t,x)=\frac{-1}{2\pi i}\int_{\Gamma (0)}e^{\lambda t}(R_{\beta
}(\lambda )f)(x)d\lambda =\frac{1}{2\pi i}\int_{\Gamma
(0)}\int_{\mathbb{R}^{3}}e^{\lambda t}\frac{e^{-\sqrt{2\lambda }|x-y|}}{2\pi |x-y|}%
g(\lambda ,y)dyd\lambda ,  \label{qq}
\end{equation}
where
\begin{equation}
g(\lambda)=(I+\beta A(\lambda ))^{-1}f. \label{qq1}
\end{equation}
By Lemma~\ref{ker10}, $A(\lambda)$ is an entire function of
$\sqrt{\lambda}$.  By the Analytic Fredholm Theorem, $(I+\beta
A(\lambda ))^{-1}$ is a meromorphic function of $\sqrt{\lambda}$,
since $A(\lambda)$ tends to zero as $\lambda \rightarrow +\infty$,
${\rm Im}(\lambda) = 0$. It does not have a pole at zero as
follows from Lemma~\ref{re1} and Remark~1 following
Lemma~\ref{ll44}. Therefore, by the Taylor formula, for all
sufficiently small $|\lambda|$, $\lambda \in \Gamma(0)$, and some
$c > 0$, we have
\begin{equation} \label{taylor}
g(\lambda)=g_{0}+g_{1}(\lambda),~~~||g_{1}(\lambda)||_{C_{\exp}(
\mathbb{R}^3)} \leq c \sqrt{|\lambda |} ||f||_{ C_{\exp}(
\mathbb{R}^3)},
\end{equation}
where $g_0 = (I+\beta A(0 ))^{-1}f$. Since $||(I+\beta A(\lambda
))^{-1}||_{C_{\exp}( \mathbb{R}^3)}$ is bounded on $\Gamma(0)$,
formula~(\ref{taylor}) is valid for all $\lambda \in \Gamma(0)$,
but not only in a neighborhood of zero.

 Let $u_{\beta }^{(1)}(x)$ be given by (\ref{qq}) with $g$
replaced by $g_{1}.$ Then
\[
u_{\beta }^{(1)}(t,x)= \frac{1}{2\pi i}\int_{\Gamma (0)}\int_{|y|
\leq \varepsilon \sqrt{t} /2}e^{\lambda t}\frac{e^{-\sqrt{2\lambda
}|x-y|}}{2\pi |x-y|} g_1(\lambda ,y)dyd\lambda +
\]
\[
\frac{1}{2\pi i}\int_{\Gamma (0)}\int_{|y|>\varepsilon \sqrt{t}
/2}e^{\lambda t}\frac{e^{-\sqrt{2\lambda }|x-y|}}{2\pi |x-y|}
g_1(\lambda ,y)dyd\lambda = I_{1}+I_{2}.
\]
We change the variable  $\lambda t=\zeta $ and use the estimate
${1}/{|x-y|}< {2}/{(\varepsilon \sqrt{t})}$ in $I_{1}$. This
implies
\[
|I_{1}|\leq \frac{c ||f||_{ C_{\exp}( \mathbb{R}^3)}}{2 \pi^2
\varepsilon t^{2}}\int_{\Gamma (0)}\int_{|y|\leq\varepsilon
\sqrt{t}/2}| \sqrt{|\zeta|} e^{\zeta -\sqrt{2\zeta }\frac{|x-y|}{
\sqrt{t}}}e^{-y^{2}} | d yd\zeta  \leq
\]
\[
 \frac{C(\varepsilon )
||f||_{ C_{\exp}( \mathbb{R}^3)} }{t^{2}},\text{ \ \ } \varepsilon
\sqrt{t} \leq |x| \leq \varepsilon ^{-1} \sqrt{t}.
\]
In $I_{2}$ we  change the variables $\lambda t=\zeta ,$
$x=\sqrt{t}z,$ \ $y= \sqrt{t}u$\  and use the estimate
$e^{-y^{2}}\leq e^{-({\varepsilon t}/{2 })^{2}}.$ This leads to
the exponential decay of $|I_{2}|$ as $t\rightarrow \infty$. Hence
\begin{equation} \label{ubet}
u_{\beta }(t,x)=\frac{1}{2\pi i}\int_{\mathbb{R}^{3}}\int_{\Gamma
(0)}e^{\lambda t} \frac{e^{-\sqrt{2\lambda }|x-y|}}{2\pi
|x-y|}g_{0}(y)d\lambda dy+r_1(t,x),
\end{equation}
where the remainder $r_1(t,x)$ satisfies
\begin{equation} \label{rdr1}
\sup_{\varepsilon \sqrt{t} \leq |x| \leq \varepsilon ^{-1}
\sqrt{t}} |r_1(t,x)| = ||f||_{ C_{\exp}( \mathbb{R}^3)}
O(t^{-2})~~{\rm as}~~t \rightarrow \infty.
\end{equation}
The integral over $\Gamma (0)$ in (\ref{ubet}) can be evaluated,
and we obtain
\[
u_{\beta }(t,x)=\frac{1}{(2\pi t)^{3/2}} \int_{\mathbb{R}^{3}}e^{-\frac{|x-y|^{2}}{2t}%
}g_{0}(y)dy+r_1(t,x).
\]
Since $||g_0||_{ C_{\exp}( \mathbb{R}^3)} \leq C  ||f||_{
C_{\exp}( \mathbb{R}^3)}$ for some constant $C$, we have
\[
u_{\beta }(t,x)= \frac{1}{(2\pi t)^{3/2}}e^{-\frac{|x|^{2}}{2t}
}\int_{\mathbb{R}^{3}}g_{0}(y)dy+r_2(t,x),
\]
where $r_2$ satisfies (\ref{rdr1}) with $r_1$ replaced by $r_2$.
In order to prove (\ref{euu0}), it remains to show that
\begin{equation} \label{qq2}
\int_{\mathbb{R}^{3}}g_{0}(x)d
x=\int_{\mathbb{R}^{3}}(1+\varphi_{\beta }(x))f(x)d x.
\end{equation}
%Denote
%\begin{equation} \label{qq3}
%w(x)=R_{0}(0)g_{0}=-\int_{\mathbb{R}^{3}}\frac{1}{2\pi
%|x-y|}g_{0}(y)d y.
%\end{equation}

Since $(I+\beta v R_0(0))g_{0}=f,$ we have $g_{0}=f-\beta v R_0(0)g_0$. Recall that $\varphi_{\beta}$
is the solution of (\ref{fibeta}) with $f=-\beta v$. Thus
\[
\int_{\mathbb{R}^{3}}g_{0}(x)d x=\int_{\mathbb{R}^{3}}f(x)dx+\int_{\mathbb{R}^{3}}
[\frac{1}{2}\Delta \varphi _{\beta }+\beta v\varphi_{\beta}]R_0(0)g_0dx.
\]
%$$g_{0}+\beta vw=f.$ Thus
%\[
%\int_{\mathbb{R}^{3}}g_{0}(x)d x=\int_{\mathbb{R}^{3}}f(x) d
%x-\int_{\mathbb{R}^{3}}\beta v w d x=\int_{\mathbb{R}^{3}}f(x)d
%x+\int_{\mathbb{R}^{3}}(\Delta \varphi _{\beta }+\beta v\varphi
%_{\beta })w d x,
%\]
%where the second equality follows from the definition of
%$\varphi_\beta$ which can be found before Theorem~\ref{thz}.
 Since $\varphi _{\beta },R_0(0)g_0=O(1/|x|)$  and their derivatives are of
  order $ O(|x|^{-2})$
as $ |x|\rightarrow \infty $, the Green formula implies
\[
\int_{\mathbb{R}^{3}}\frac{1}{2}\Delta \varphi _{\beta }R_0(0)g_0 dx=
\int_{\mathbb{R}^{3}}\varphi _{\beta }\frac{1}{2}\Delta R_0(0)g_0 d x=
\int_{\mathbb{R}^{3}}\varphi _{\beta }g_0 d x.
\]
Hence
%\begin{equation} \label{qq4}
\[
\int_{\mathbb{R}^{3}}g_{0}(x)d x=\int_{\mathbb{R}^{3}}f(x)d
x+\int_{\mathbb{R}^{3}}\varphi _{\beta }(I+\beta v R_0(0))g_{0}d x,
\]
which implies (\ref{qq2}.)
%\end{equation}
%From (\ref{aa3}) and (\ref{qq3}) it follows that
%\[
%(\Delta +\beta v)w=(\Delta +\beta v)[R_0(0)(1+\beta
%A(0))^{-1}f]=f.
%\]
%Thus (\ref{qq4}) implies (\ref{qq2}).
 \qed

 Next, let us study the distribution of the polymer with
respect to the measure $\mathrm{P}_{\beta,T}$ as $T \rightarrow
\infty$. Consider the process $y^T(t) = x(t T)/\sqrt{T}$, $0 \leq
t \leq 1$.
\begin{theorem} \label{nn1} Let $d \geq 3$ and  $0 \leq \beta <
\beta_{cr}$. With respect to $\mathrm{P}_{\beta,T}$, the
distribution of the process $y^T(t)$ converges as $T \rightarrow
\infty$, weakly in the space $C([0,1], \mathbb{R}^d)$, to the
distribution of the $d$-dimensional Brownian motion. With respect
to $\mathrm{P}_{\beta,T}(\cdot|x(T)=0)$, the distribution of the
process $y^T(t)$ converges as $T \rightarrow \infty$, weakly in
the space $C([0,1], \mathbb{R}^d)$, to the distribution of the
$d$-dimensional Brownian bridge.
\end{theorem}
\proof We shall only prove the first statement since the proof of
the second one is completely similar. First, let us prove the
convergence of the finite-dimensional distributions. Clearly
$\mathrm{P}_{\beta,T}(y^T(0) = 0) =1$.
 Let $0 < t_1 < ... < t_n \leq 1$. The density of the random
vector $(y^T(t_1),...,y^T(t_n))$ with respect to the Lebesgue
measure on $ \mathbb{R}^{dn}$ is equal to
\[
\rho^T(x_1,...,x_n) =
\]
\[
 T^{\frac{dn}{2}}p_\beta(t_1 T,  0, x_1 T^{\frac{1}{2}})
  p_\beta((t_2- t_1)T, x_1T^{\frac{1}{2}}, x_2T^{\frac{1}{2}})...
p_\beta((t_n - t_{n-1})T,x_{n-1}T^{\frac{1}{2}},
x_{n}T^{\frac{1}{2}}) {(Z_{\beta, T}(0))^{-1}}.
\]
By Lemma~\ref{uu1},
\[
p_\beta(t_1 T,  0, x_1 T^{\frac{1}{2}}) = T^{-d/2} (2 \pi
t_1)^{-d/2} (1 + \varphi_\beta(0)) \exp(-|x_1|^2/2t_1)(1 +
r(T,x_1)),
\]
where
\begin{equation} \label{rx1}
 \lim_{T \rightarrow \infty} \sup_{\varepsilon  \leq |x_1|
\leq \varepsilon^{-1}} (|r(T,x_1|) = 0.
\end{equation}
Note that $p_\beta \geq p_0$ since $v$ is non-negative, and
$\lim_{T \rightarrow \infty} (Z_{\beta, T}(0)) = (1 +
\varphi_\beta(0))$ by Theorem~\ref{thz}. Therefore,
\[
\rho^T(x_1,...,x_n) \geq
\]
\begin{equation} \label{osi}
(2 \pi t_1)^{-\frac{d}{2}} e^{-\frac{|x_1|^2}{2t_1}}(1 + r(T,x_1))
(2 \pi (t_2-t_1))^{-\frac{d}{2}}
e^{-\frac{|x_2-x_1|^2}{2(t_2-t_1)}}... (2 \pi
(t_n-t_{n-1}))^{-\frac{d}{2}}
e^{-\frac{|x_n-x_{n-1}|^2}{2(t_n-t_{n-1})}}
\end{equation}
\[
=\rho^W_{t_1,...,t_n}(x_1,...,x_n)(1 + r(T,x_1)),
\]
where ${\rho}^W_{t_1,...,t_n}(x_1,...,x_n)$ is the density of the
Gaussian vector $(W(t_1),...,W(t_n))$, where $W$ is a
$d$-dimensional Brownian motion,  and $q(T,x_1)$ satisfies
(\ref{rx1}) with $q$ instead of $r$. Since $\varepsilon$ was an
arbitrary positive number, this implies the convergence of the
finite-dimensional distributions of $y^T$ to the
finite-dimensional distributions  of the Brownian motion. Indeed,
the estimate from below for $\rho^T(x_1,...,x_n)$ in (\ref{osi})
is sufficient since we know a priori that
$\rho^W_{t_1,...,t_n}(x_1,...,x_n)$ is the density of a
probability measure.

It remains to prove tightness of the family of processes $y^T$, $T
\geq 1$.

For a
 continuous function $x: [0,T]
\rightarrow {\mathbb{R}}^d$, let
\[
m(x, \delta) = \sup_{|t_1  - t_2| \leq \delta T,~ 0 \leq t_1, t_2
\leq T} |x(t_1) - x(t_2)|/\sqrt{T},
\]
\[
\widetilde{m}(x, \delta, \varepsilon) = \sup_{|t_1 - t_2| \leq
\delta T,~ 0 \leq t_1, t_2 \leq T,~|x(t_1)| \geq \varepsilon
\sqrt{T}} |x(t_1) - x(t_2)|/\sqrt{T}.
\]
The tightness will follow if we show that for each $\varepsilon,
\eta > 0$ there is $\delta > 0$ such that
\[
\mathrm{P}_{\beta,T}(m(x, \delta) > \varepsilon) \leq \eta
\]
for all sufficiently large $T$. Note that $m(x, \delta) >
\varepsilon$ implies that $\widetilde{m}(x, \delta, \varepsilon/4)
> \varepsilon/4$. Therefore, it is sufficient to show that
\begin{equation} \label{tight01}
\mathrm{P}_{\beta,T}(\widetilde{m}(x, \delta, \varepsilon/4) >
\varepsilon/4) \leq \eta.
\end{equation}
Fix $\varepsilon > 0$.
%Since the support of $v$ is compact, we can
%find $R > 0$ such that the support of $v$ is contained in the ball
%of radius $R$ centered at the origin.
For a continuous function
$x: [0,T] \rightarrow \mathbb{R}^d$, let
\[
\tau = \min(T, \inf\{t \geq 0: |x(t)| = \varepsilon \sqrt{T}/4\}
),
\]
 Let $\mathcal{E}_\delta$ be the event that
${m}(x, \delta) > \varepsilon/4$ and
  $ \widetilde{\mathcal{E}}_\delta$ the event that
$\widetilde{m}(x, \delta, \varepsilon/4) > \varepsilon/4$. For $0
\leq s \leq T$, let  $\mathcal{E}^s_\delta$ be the event that a
continuous function $x: [0,T-s] \rightarrow {\mathbb{R}}^d$
satisfies
\[
\sup_{|t_1  - t_2| \leq \delta T,~ 0 \leq t_1, t_2 \leq T -s}
|x(t_1) - x(t_2)|/\sqrt{T} > \varepsilon/4.
\]
Then
\[
\mathrm{P}_{\beta,T}( \widetilde{\mathcal{E}}_\delta) = (Z_{\beta,
T}(0))^{-1} \mathrm{E}_{0,T}(\exp(\int_0^{T} \beta v (x(t)) d t)
\chi_{ \widetilde{\mathcal{E}}_\delta}) \leq
\]
\[
 (Z_{\beta, T}(0))^{-1}
\mathrm{E}_{0,T}\left(\exp(\int_0^{\tau} \beta v (x(t)) d t)
\mathrm{E}^{x(\tau)}_{0,T-\tau} ( \chi_{ \mathcal{E}^\tau_\delta}
\exp( \int_0^{T -\tau} \beta v (x(t)) d t ) ) \right),
\]
where $ \mathrm{E}^{x}_{0,T}$ denotes the expectation with respect
to the measure induced by the Brownian motion starting at the
point $x$.
 Since
\[
\mathrm{E}_{0,T} \exp(\int_0^{\tau} \beta v (x(t)) d t) \leq
Z_{\beta, T}(0)
\]
and
\[
\mathrm{E}^{x(\tau)}_{0,T-\tau} ( \chi_{ \mathcal{E}^\tau_\delta}
\exp( \int_0^{T -\tau} \beta v (x(t)) d t ) ) \leq \sup_{x \in
\mathbb{R}^d, |x| = \varepsilon \sqrt{T}/4} \mathrm{E}^{x}_{0,T} (
\chi_{ \mathcal{E}_\delta} \exp( \int_0^{T} \beta v (x(t)) d t )
),
\]
 it is sufficient to estimate
\begin{equation} \label{sup22}
\sup_{x \in \mathbb{R}^d, |x| = \varepsilon \sqrt{T}/4}
\mathrm{E}^{x}_{0,T} ( \chi_{ \mathcal{E}_\delta} \exp (\int_0^{T}
\beta v (x(t)) d t )).
\end{equation}
Let $ \mathcal{E}'$ be the event that a trajectory starting at $x$
reaches the support of $v$ before time $T$. Note that
\[
 \lim_{T
\rightarrow \infty}  \sup_{x \in \mathbb{R}^d, |x| = \varepsilon
\sqrt{T}/4} \mathrm{P}^x_{0,T} (\mathcal{E}') = 0
\]
 since $d \geq
3$. The expression in (\ref{sup22}) is estimated form above by
\[
\sup_{x \in \mathbb{R}^d, |x| = \varepsilon \sqrt{T}/4} (
\mathrm{E}^{x}_{0,T} ( \chi_{ \mathcal{E}'} \exp( \int_0^{T} \beta
v (x(t)) d t )) + \mathrm{P}^{x}_{0,T} ({ \mathcal{E}_\delta})).
\]
The second term does not depend on $T$ due to the scaling
invariance of the Brownian motion, and can be made arbitrarily
small by selecting a sufficiently small $\delta$. Due to the
Markov property of the Brownian motion, the first term is
estimated from above by
\[
 \sup_{x
\in \mathbb{R}^d, |x| = \varepsilon T/4} \mathrm{P}^x_{0,T}
(\mathcal{E}') \cdot \sup_{x \in {\rm supp} (v)} Z_{\beta, T}(x),
\]
and thus tends to zero when $T \rightarrow \infty$. \qed

\section{Behavior of the Polymer for $\protect\beta = \protect%
\beta_{cr}$} \label{bb2}

In this section we assume that $d = 3$. Again, we start with the
asymptotic behavior of the solution $u_\beta(t,x)$ of the Cauchy
problem and of the fundamental solution $p_\beta(t,y,x)$ when $t
\rightarrow \infty$, $|y| \leq \varepsilon^{-1}$, $\varepsilon
\sqrt{t} \leq |x| \leq \varepsilon^{-1} \sqrt{t}$, and
$\varepsilon > 0$ is small but fixed.

Recall that $\psi$ is the positive ground state for
$H_{\beta_{cr}}$ normalized by the condition ${||\beta_{cr} v
\psi||_{L^2_{\exp}(\mathbb{R}^3)} = 1}$ (see the remark following
Lemma~\ref{re1} and Theorem~\ref{thz}). For $f \in C_{\exp}(
\mathbb{R}^3)$, define
\[
\alpha(f) = \varkappa \int_{\mathbb{R}^3}\psi (x)f(x)d x, ~~~
\varkappa = \frac{ 1} {
\sqrt {2\pi}  \beta_{cr} \int _{\mathbb{R}^3}v(x)\psi (x)d
x}.
\]
We can formally apply
this to $f$ being the $\delta$-function centered at a point $y$,
and thus define
\[
\alpha(\delta_y(x)) = \varkappa \psi (y).
\]

\begin{theorem} \label{uu1xx}
Let $d = 3$,  $ \beta = \beta_{cr}$, $\varepsilon > 0$ and $f \in
C_{\exp}( \mathbb{R}^3)$, $f \geq 0$. We have the following
asymptotics for the solution $u_{\beta}$ of the Cauchy problem
with the initial data $f$:
\begin{equation} \label{euu0xx}
u_{\beta}(t,x) = \frac{1}{|x| \sqrt{t}} \exp(-|x|^2/2t) (
\alpha(f) + q_f(t,x)) ,
\end{equation}
where for some constant $C_\beta(\varepsilon)$ we have
\[
 \sup_{\varepsilon \sqrt{t} \leq |x|
\leq \varepsilon^{-1} \sqrt{t}}| q_f(t,x)| \leq
C_\beta(\varepsilon) t^{-1/2} ||f||_{ C_{\exp}( \mathbb{R}^3)},~~t
\geq 1.
\]

We have the following asymptotics for the fundamental solution of
the parabolic equation:
\begin{equation} \label{euu1xx}
p_\beta(t,y,x) =  \frac{\varkappa}{|x| \sqrt{t}} \exp(-|x|^2/2t)(
\psi(y) + q(t,y,x)),
\end{equation}
where
\[
 \lim_{t \rightarrow \infty} \sup_{|y| \leq \varepsilon^{-1},~
\varepsilon \sqrt{t} \leq |x| \leq \varepsilon^{-1} \sqrt{t}}
|q(t,y,x| = 0.
\]
\end{theorem}
\proof As in Lemma~\ref{uu1}, formula (\ref{euu1xx}) follows from
(\ref{euu0xx}). Lemma \ref{tet} implies
\[
(I + \beta_{cr} A(\lambda))^{-1} = \frac{\beta_{cr}
}{\beta(\lambda) - \beta_{cr} }B + O(1),~~~\lambda \rightarrow
0,~\lambda \in \mathbb{C}',
\]
where $B$ is the one dimensional operator with the kernel
\[
B(x,y) = \frac{v(x) \psi (x) \psi (y)}{ \int
_{\mathbb{R}^3}v(x)\psi ^2(x)dx}.
\]
From here,  (\ref{dim3}) and (\ref{betas}) we get
\[
(I + \beta_{cr} A(\lambda))^{-1} = \frac{1}{\beta_{cr} \gamma
\sqrt \lambda }B + O(1),~~~\lambda \rightarrow 0,~\lambda \in
\mathbb{C}',
\]
where $\gamma$ is defined in (\ref{gamma}), (\ref{gamma3}). Hence, for any $f \in
C_{\exp}( \mathbb{R}^3)$ and $\lambda \rightarrow 0$, $\lambda \in
\mathbb{C}'$,
\begin{equation} \label{h00}
h(\lambda, x) :=(I + \beta_{cr} A(\lambda))^{-1}f =
\frac{\widetilde{\alpha}(f)}{\sqrt \lambda } v(x) \psi (x)+
g_1(\lambda),~~~\widetilde{\alpha}(f)=\frac{\sqrt {2}\pi
\int_{\mathbb{R}^3}\psi (x)f(x)d x}{\beta_{cr}
(\int_{\mathbb{R}^3}v(x)\psi (x)d x)^2},
\end{equation}
where $g_1(\lambda) \leq c ||f||_{C_{\exp}( \mathbb{R}^3)}$ for
some constant $c$. Now, similarly to (\ref{qq}), we have
\[
u_{\beta }(t,x)=\frac{-1}{2\pi i}\int_{\Gamma (0)}e^{\lambda t}(R_{\beta
}(\lambda )f)(x)d\lambda =\frac{1}{2\pi i}\int_{\Gamma
(0)}\int_{\mathbb{R}^{3}}e^{\lambda t}\frac{e^{-\sqrt{2\lambda }|x-y|}}{2\pi |x-y|}%
h(\lambda ,y)dyd\lambda.
\]
The integral with $g_1(\lambda)$ instead of $h$ can be estimated
similarly to the estimate on $u_{\beta}^{(1)}$ in the case of
$\beta < \beta_{cr}$. This leads to following analogue of
(\ref{ubet})
\[
u_{\beta }(t,x)=\frac{\widetilde{\alpha }(f)}{2\pi
i}\int_{\mathbb{R}^{3}}\int_{\Gamma (0)}e^{\lambda t}
\frac{e^{-\sqrt{2\lambda}|x-y|}}{2\pi \sqrt{\lambda} |x-y|}v(y)
\psi (y)d\lambda d y+r_1(t,x),
\]
where the remainder $r_1(t,x)$ satisfies
\begin{equation} \label{rdr1xx}
\sup_{\varepsilon \sqrt{t} \leq |x| \leq \varepsilon ^{-1}
\sqrt{t}} |r_1(t,x)| = ||f||_{ C_{\exp}( \mathbb{R}^3)}
O(t^{-3/2})~~{\rm as}~~t \rightarrow \infty.
\end{equation}
We evaluate the integral over $\Gamma (0)$:
\begin{equation} \label{lapt}
\frac{1}{2\pi i}\int_{\Gamma (0)} \frac{e^{\lambda t -
\sqrt{2\lambda }|x-y|}}{ \sqrt{2\lambda} }d \lambda =   \frac
{1}{\sqrt {2\pi t}}e^{-\frac{|x-y|^2}{2t}}.
\end{equation}
This equality simply means that the  inverse Laplace transform of the Green function of the one
dimensional Helmholtz equation coincides
with the fundamental solution of the corresponding heat equation. Thus,
\[
u_{\beta }(t,x)=\frac{\widetilde{\alpha}(f)}{2\pi^{3/2} \sqrt {t}}
\int_{\mathbb{R}^3} \frac {1}{|x-y|} e^{-\frac {|x-y|^2}{2t}}v(y)
\psi(y)dy+r_1(t,x).
\]
This implies (\ref{euu0xx}) since $v$ has a compact support. \qed

The next theorem concerns the fundamental solution when both $y$
and $x$ are at a distance of order $\sqrt{t}$ away from the
origin. Note that now there are two terms in the asymptotic
expansion for the fundamental solution which are of the same order
in $t$. The main terms have the order $t^{-3/2}$ when $t
\rightarrow \infty$, compared with $t^{-1}$ in the case considered
in Theorem~\ref{uu1xx} (where $y$ was bounded).

%Denote by $c_{\beta}$ the coefficient in the asymptotic behavior of the ground state at infinity:
%\[
%c_{\beta}=\lim_{|x|\rightarrow\infty}(|x|\psi(x)).
%\]
%The existence and positivity of the limit easily follow from the
%proof of Lemma~\ref{re1}.

\begin{theorem} \label{uu1xx1}
Let $d = 3$,  $ \beta = \beta_{cr}$, $\varepsilon > 0$.
We have the following asymptotics for the fundamental solution of
the parabolic equation:
\begin{equation} \label{euu1xx1}
p_\beta(t,y,x) = p_0(t,y,x)+ \frac{1}{(2\pi)^{3/2}|y||x| \sqrt{t}}
e^{-(|y|+|x|)^2/2t}( 1 + \overline{q}(t,y,x)),
\end{equation}
where
\begin{equation} \label{qbar}
 \lim_{t \rightarrow \infty} \sup_{
\varepsilon \sqrt{t} \leq |y|,|x| \leq \varepsilon^{-1} \sqrt{t}}
|\overline{q}(t,y,x| = 0.
\end{equation}
\end{theorem}
\proof Let $p_\beta(t,y,x) = p_0(t,y,x)+u$. Then
$u_t=H_{\beta}u+\beta vp_0,~~u|_{t=0}=0$, and therefore by the
Duhamel formula
\[
u(t,y,x)=\int_0^t \int_{\mathbb{R}^3}p_\beta(t-s,z,x)\beta
v(z)p_0(s,y,z)dzds.
\]
Using (\ref{euu1xx}), we get
\begin{equation} \label{int11a}
u(t,y,x )=\int_0^t \int_{\mathbb{R}^3}\frac{\varkappa}{|x|
\sqrt{t-s}} \exp(-|x|^2/2(t-s))( \psi(z)\beta v(z)p_0(s,y,0)dzds +
h_1 + h_2
\end{equation}
with
\[
h_1 = \int_0^t \int_{\mathbb{R}^3}\frac{\varkappa}{|x| \sqrt{t-s}}
\exp(-|x|^2/2(t-s))( \psi(z)\beta v(z)(p_0(s,y,z)-
p_0(s,y,0))dzds,
\]
\[
h_2=\int_0^t \int_{\mathbb{R}^3}\frac{\varkappa}{|x| \sqrt{t-s}}
\exp(-|x|^2/2(t-s))( q(t-s,z,x)\beta v(z)p_0(s,y,z)dzds,
\]
where $q$ is the same as in (\ref{euu1xx}). The integral in the
right hand side of (\ref{int11a}) (let us denote it by $w$) is a
convolution of two functions and can be evaluated using the
Laplace transform (see (\ref{lapt})). It gives the second term in
the right hand side of (\ref{euu1xx1}). The contribution from the
other two terms can be shown to satisfy~(\ref{qbar}). Let us prove
the statement about $w$. In fact,
\[
w=\varkappa_1(w_1*p_0(t,y,0)),~~ \varkappa_1=\frac{\varkappa \sqrt
{2\pi}}{|x|}\int_{\mathbb{R}^3} \beta v(z)
\psi(z)dz=\frac{1}{|x|},~~w_1= \frac{1} {\sqrt {2\pi t}}
\exp(-|x|^2/2t)).
\]
The Laplace transform $\widehat{w}_1(\lambda)$ of the function
$w_1$ is equal to $e^{-\sqrt{2\lambda} |x|}/\sqrt{2\lambda}$ (see
(\ref{lapt})), and the Laplace transform of $p_0(t,y,0)$ is equal
to $e^{-\sqrt{2\lambda} |y|}/2\pi |y|$. Thus
\[
\widehat{w}(\lambda)=\frac{1}{2\pi |x||y|}\frac{
e^{-\sqrt{2\lambda} (|x|+|y|)}}{ \sqrt{2\lambda}}.
\]
It remains to apply (\ref{lapt}) one more time.
 \qed
\\

As in Section~\ref{bb1x}, we shall study the limit, as $T
\rightarrow \infty$, of the family of processes
 $y^T(t) = x(t T)/\sqrt{T}$, $0 \leq
t \leq 1$.
 For $0 \leq s < t \leq 1$, $y,x \in \mathbb{R}^3$, define
\[
p_\beta^T(s,t,y,x) = p_\beta(T (t-s), y \sqrt{T}, x \sqrt{T} ),
\]
\[
\overline{p}_\beta(s,t,0,x) =\lim_{T \rightarrow \infty} (T
{p}^T_\beta(s,t,0,x)) =  \lim_{T \rightarrow \infty} (T p_\beta(T
(t-s), 0, x \sqrt{T} )),~x \neq 0,
\]
\[
{\overline{p}}_\beta(s,t,y,x) = \lim_{T \rightarrow \infty}
(T^{3/2} {{p}}^T_\beta(s,t,y,x)) =  \lim_{T \rightarrow \infty} (
T^{3/2} p_\beta(T(t-s), y \sqrt{T}, x \sqrt{T} )),~~y,x \neq 0.
\]
By Theorems~\ref{uu1xx} and \ref{uu1xx1},
\[
\overline{p}_\beta(s,t,0,x) = \frac{\varkappa \psi(0)}{|x|
\sqrt{t-s}} \exp(-|x|^2/2(t-s)),~~x \neq 0,
\]
\begin{equation} \label{pch}
{\overline{p}}_\beta(s,t,y,x) =  p_0(t-s,y,x)+ \frac{1}
{(2\pi)^{3/2}|y||x| \sqrt{t-s}} \exp(-(|y|+|x|)^2/2(t-s))~~,y,x \neq 0.
\end{equation}
For $0 < t_1 < ... < t_n \leq 1$, let the density of the random
vector $(y^T(t_1),...,y^T(t_n))$ with respect to the Lebesgue
measure on $ \mathbb{R}^{dn}$ be denoted by $
\rho^T(x_1,...,x_n)$.

For $0 \leq s < t \leq 1$ and $y, x \in \mathbb{R}^3 $, define
\[
Q^T(s,t,y,x) = {p}^T_\beta(s,t,y,x) \int_{ \mathbb{R}^3}
{{p}}^T_\beta(t,1,x,z) d z (\int_{ \mathbb{R}^3}
{p}^T_\beta(s,1,y,z) d z )^{-1},~~t <1,
\]
\[
Q^T(s,1,y,x) = {p}^T_\beta(s,1,y,x) (\int_{ \mathbb{R}^3}
{p}^T_\beta(s,1,y,z) d z )^{-1}.
\]
Thus
\[
\rho^T(x_1,...,x_n) =
Q^T(0,t_1,0,x_1)Q^T(t_1,t_2,x_1,x_2)...Q^T(t_{n-1}, t_n, x_{n-1},
x_n).
\]
In order to find the limit of the finite dimensional distributions
of $y^T$, we need to identify the limit of $Q^T$ as $T \rightarrow
\infty$. For $0 \leq s < t \leq 1$, $y \in \mathbb{R}^3$ and $ x
\in \mathbb{R}^3 \setminus \{ 0\} $, define
\begin{equation} \label{qdefe}
Q(s,t,y,x) = \lim_{T \rightarrow \infty} Q^T(s,t,y,x).
\end{equation}
 By
Theorems~\ref{uu1xx} and \ref{uu1xx1},
\begin{equation} \label{q11}
 Q(s,t,y,x)
={\overline{p}}_\beta(s,t,y,x) \int_{ \mathbb{R}^3}
{\overline{p}}_\beta(t,1,x,z) d z (\int_{ \mathbb{R}^3}
{\overline{p}}_\beta(s,1,y,z) d z )^{-1},~~t <1,
\end{equation}
\begin{equation} \label{q12}
Q(s,1,y,x) = {\overline{p}}_\beta(s,1,y,x) (\int_{ \mathbb{R}^3}
{\overline{p}}_\beta(s,1,y,z) d z )^{-1}.
\end{equation}
 We additionally define
$Q(s,t,y,0) = 0$.

Using (\ref{qdefe}), (\ref{q11}) and (\ref{q12}), we can identify
the limit of the densities $\rho^T(x_1,...,x_n)$ for $x_2,...,x_n
\neq 0$. In order to identify the weak limit of the finite
dimensional distributions of the processes $y^T$, we are going to
show that the limit of the densities is the density of a
probability distribution, i.e. the mass does not escape to the
origin or infinity. This is done in Lemma~\ref{qpr}, where we show
that $Q$ serves as the transition density for a Markov process.
First, however, we show that $Q$ satisfies a Fokker-Plank type
equation on $ \mathbb{R}^3 \setminus \{0 \}$.
% with the operator in
%the right hand side being the formal adjoint of the generator of
%the Markov process.

Let
\begin{equation} \label{gtx}
g(t,x) = \ln(\int_{ \mathbb{R}^3} \overline{p}_\beta(t,1,x,z) d
z),~~0 \leq t < 1,~~|x| > 0,
\end{equation}
Let $L$ be the differential operator acting on $C^2( \mathbb{R}^3
\setminus \{0\})$ according to the formula
\[
(L f) (t,x) = \frac{1}{2} \Delta_x f (t,x) + (\frac{\partial
g(t,x)}{\partial r}) \frac{\partial f}{\partial r} (t,x),~~~|x| >
0,
\]
and let $L^*$ be the formal adjoint of $L$, i.e.
\[
L^*v=\frac{1}{2} \Delta_x v- \frac{1}{r^2}
\frac{\partial[(\partial g /\partial r) v]}{ \partial r}.
\]
\begin{lemma}  For $0 \leq s < 1$ and $y \in \mathbb{R}^3$, the function $
Q(s,t,y,x)$ satisfies the equation
\begin{equation} \label{L*}
\frac{\partial Q(s,t,y,x) }{\partial t} = L^* Q(s,t,y,x),~~~|x| >
0,~~s < t < 1.
\end{equation}
\end{lemma}
\proof Let us consider the case when $y \neq 0$ (the other case is
similar). Let
\[
v_1(s,t,y,x) = \frac{1}{(2\pi)^{3/2}|y||x| \sqrt{t-s}}
\exp(-(|y|+|x|)^2/2(t-s)),
\]
\[
v_2(t,x) = \int_{ \mathbb{R}^3} \frac{1}{(2\pi)^{3/2}|x||z|
\sqrt{1-t}} \exp(-(|x|+|z|)^2/2(1-t)) d z.
\]
%
%\begin{equation*}
%p_{\beta }(t,y,x)=p_{0}(t,y,x)+v_{1}(t,y,x),
%\end{equation*}
Observe that
\begin{equation} \label{diff33}
(\frac{\partial }{\partial t}-\frac{1}{2}\Delta_x
)v_{1}=0,~~~(\frac{\partial }{\partial t}+\frac{1}{2}\Delta_x
)v_{2}=0.
\end{equation}
For fixed $s$ and $y$, the function $Q(s,t,y,x)$ is proportional
to
%The limiting \ density distribution is proportional to
\[
u(t,x)=
%(p_{0}(t,y,x)+v_{1}(t,y,x))[1+\int
%v_{1}(1-t,x,z)dz]=
(p_{0}(t-s,y,x)+v_{1}(s,t,y,x))[1+v_{2}(t,x)].
\]
By (\ref{diff33}),
\begin{equation} \label{hhw1}
(\frac{\partial }{\partial t}-\frac{1}{2}\Delta_x
)u=-(\frac{\partial p_{0}}{\partial r} + \frac{\partial
v_{1}}{\partial r}) \frac{\partial v_{2}}{\partial r}
+2(p_{0}+v_{1}) \frac{\partial v_{2}}{\partial t}.
\end{equation}
%Let
%\[
%A=\frac{\partial v_{2}}{\partial r} (1+v_{2})^{-1},~~~ B= -(2
%\frac{\partial v_{2}}{\partial t} +A \frac{\partial
%v_{2}}{\partial r})( 1+v_{2})^{-1}.
%\]
%If we apply  the operator formally  adjoint  to $ (\frac{\partial
%g(t,x)}{\partial r}) \frac{\partial }{\partial r}$ to the function
%$u$, we obtain
For any two functions $A$ and $B$ we have
\begin{equation} \label{hhw2}
(A\frac{\partial }{\partial r} +B)u=A(\frac{\partial
p_{0}}{\partial r} +\frac{\partial v_{1}}{\partial r}
)(1+v_{2})+A(p_{0}+v_{1})\frac{\partial v_{2}}{\partial r}
+B(p_{0}+v_{1})(1+v_{2}).
\end{equation}
Thus
\[
(\frac{\partial }{\partial t}-\frac{1}{2}\Delta_x +A\frac{\partial
}{\partial r} +B)u=
\]
\[
(\frac{\partial p_{0}}{\partial r} + \frac{\partial
v_{1}}{\partial r}) (-\frac{\partial v_{2}}{\partial r}+A(1+v_2))
+2(p_{0}+v_{1}) (\frac{\partial v_{2}}{\partial t}+A\frac{\partial
v_{2}}{\partial r} +B(1+v_{2}))=0
\]
if
\[
A=\frac{\partial v_{2}}{\partial r} (1+v_{2})^{-1},~~~ B= -(2
\frac{\partial v_{2}}{\partial t} +A \frac{\partial
v_{2}}{\partial r})( 1+v_{2})^{-1}.
\]
Since $g(t,x) = \ln (1+v_2)$ and $2 {\partial v_{2}}/{\partial
t}=- {\partial^2 v_{2}}/{\partial r^2}-2{\partial v_{2}}/{\partial
r}$ (see (\ref{diff33})), it is easy to check that the operator in
the left hand side of the equation for $u$ is $\frac{\partial
}{\partial t}-L^* $, and this justifies (\ref{L*}). \qed
%Adding (\ref{hhw1}) and (\ref{hhw2}), we obtain
%\[
%\frac{\partial u(t,x) }{\partial t} = L^* u(t,x),~~~|x| > 0,
%\]
%as claimed.
\begin{lemma} \label{qpr} The function $Q(s,t,y,x)$, $0 \leq s < t \leq 1$,
$y,x \in \mathbb{R}^3$, is the transition density for a Markov
process on $ \mathbb{R}^3$.
\end{lemma}
\proof To show the existence of a Markov process, we need to
verify that
\begin{equation} \label{cons}
\int_{ \mathbb{R}^3} Q(t_1, t_2, x_1, x_2) d x_2  = 1,~~t_1 < t_2
\end{equation}
and
\begin{equation} \label{semig}
\int_{ \mathbb{R}^3} Q(t_1, t_2, x_1, x_2) Q(t_2, t_3, x_2, x_3)
dx_2 = Q(t_1, t_3, x_1, x_3),~~t_1 < t_2 < t_3.
\end{equation}
Let us assume that (\ref{cons}) has been demonstrated, and prove
(\ref{semig}). Observe that
\[
 \int_{ \mathbb{R}^3} T^{2 + \alpha} p^T_\beta(t_1, t_2, x_1, x_2)
p^T_\beta(t_2, t_3, x_2, x_3) dx_2 = T^{1+ \alpha} p^T_\beta(t_1,
t_3, x_1, x_3),~~t_1 < t_2 < t_3,
\]
where $\alpha = 1/2$ if $x_1 = 0$ and $\alpha = 0$ otherwise.  For
$x_3 \neq 0$ we take the limit, as $T \rightarrow \infty$, on both
sides of this relation. The integrand on the left hand side
converges to
\[
\overline{p}_\beta(t_1, t_2, x_1, x_2) \overline{p}_\beta(t_2,
t_3, x_2, x_3),
\]
however, the convergence is not necessarily uniform in $x_2$, and
we can only conclude by the Fatou Lemma that
\[
 \int_{ \mathbb{R}^3 \setminus \{ 0\} } \overline{p}_\beta(t_1, t_2, x_1, x_2) \overline{p}_\beta(t_2,
t_3, x_2, x_3) dx_2 \leq  \overline{p}_\beta(t_1, t_3, x_1,
x_3),~~t_1 < t_2 < t_3,~~ x_3 \neq 0.
\]
From (\ref{q11}) and (\ref{q12}) it now follows that
\[
 \int_{ \mathbb{R}^3 \setminus \{ 0\} } Q(t_1, t_2, x_1, x_2) Q(t_2,
t_3, x_2, x_3) dx_2 \leq  Q(t_1, t_3, x_1, x_3),~~t_1 < t_2 <
t_3,~~ x_3 \neq 0.
\]
Note that both sides of this inequality are continuous in $x_3 \in
\mathbb{R}^3 \setminus \{ 0\}$.  Due to (\ref{cons}), the
integrals in $x_3$ over $\mathbb{R}^3 \setminus \{ 0\}$ are equal
to one for the expressions in both sides of this inequality.
Therefore,
\[
 \int_{ \mathbb{R}^3 \setminus \{ 0\} } Q(t_1, t_2, x_1, x_2) Q(t_2,
t_3, x_2, x_3) dx_2 =  Q(t_1, t_3, x_1, x_3),~~t_1 < t_2 < t_3,~~
x_3 \neq 0,
\]
and thus (\ref{cons}) implies (\ref{semig}).

Now let us verify (\ref{cons}). Put $s = t_1, \tau = t_2, y = x_1$
and $x = x_2$. Again, we shall consider the case $y \neq 0$, the
other case being similar. Moreover, we can assume that $\tau < 1$,
since the case $\tau = 1$ can be treated by taking the limit $\tau
\uparrow 1$.
%Add together. Terms with $p_{0_{r}}$ have to canceled and terms with $%
%p_{0_{r}}:$
%
%\begin{equation*}
%A=\frac{v_{2_{r}}}{1+v_{2}},\text{ \ \ }B=\frac{-2v_{2_{t}}-Av_{2_{r}}}{%
%1+v_{2}}.
%\end{equation*}
%We have $2v_{2_{t}}=-v_{2_{rr}}-\frac{2}{r}v_{2_{r}}.$%
%\begin{equation*}
%\frac{v_{2_{rr}}+\frac{2}{r}v_{2_{r}}-Av_{2_{r}}}{1+v_{2}}=?A^{\prime }+2A/r=%
%\frac{v_{2_{rr}}}{1+v_{2}}-(\frac{v_{2_{r}}}{1+v_{2}})^{2}+2\frac{v_{2_{r}}}{%
%1+v_{2}}/r
%\end{equation*}
On a formal level, (\ref{cons}) follows from (\ref{L*}) by
integrating the both sides of (\ref{L*}) over $\Omega
=[s,\tau]\times \mathbb{R}^3 \subset \mathbb{R}^4_{t,x}$:
\begin{equation} \label{q112}
\int_{\mathbb{R}^3}Q(s,\tau,y,x)dx- \lim_{t \downarrow s}
\int_{\mathbb{R}^3}Q(s,t,y,x)dx =\langle
L^*Q,1\rangle_{L^2(\Omega)}=\langle Q,L1\rangle_{L^2(\Omega)}.
\end{equation}
One needs only to note that
\begin{equation} \label{q111}
\lim_{t \downarrow s} \int_{\mathbb{R}^3}Q(s,t,y,x)dx=1,
\end{equation}
and that the operator $L$ applied to the identity function gives
zero. The latter implies that the left hand side in (\ref{q112})
is zero, and  (\ref{q111}) implies that the second term on the
left hand side of (\ref{q112}) is one.

In order to make relations (\ref{q112}) rigorous we note that
$Q(s,t,y,x)$ is infinitely smooth in $(t,x)$ when $x\neq0$ and
decays exponentially as $|x|\rightarrow \infty$.
 However, it has a singularity at $x=0$. Thus the integrals
over $\mathbb{R}^3$ and $\Omega$ in (\ref{q112}) must be
understood as limits of the corresponding integrals over the
region $|x|>\varepsilon$ as $\varepsilon \rightarrow 0$. Let us
examine the singularities of $Q$ and of the coefficients of $L^*$
at the origin.

Relation (\ref{pch}) implies that
\begin{equation} \label{q22}
\overline{p}_\beta(s,t, y, x)=\frac{a}{r}+O(r),~~r =
|x|\rightarrow0,~~a=a(s,t,y).
\end{equation}
It is important that (\ref{q22}) does not contain a term of order
$O(1)$. From (\ref{q22}), (\ref{gtx}) and (\ref{q11}) it follows
that
\begin{equation}
\frac{\partial g(t,x)}{r}=-\frac{1}{r}+O(r),~~ \label{last2}
Q(s,t,y,x)=\frac{c}{r^2}+O(1), ~~\frac{\partial
Q(s,t,y,x)}{\partial r} =-\frac{2c}{r^3}+O(1),
\end{equation}
where $r\rightarrow0$, $c=c(s,t,y)$. Since $Q$ has a weak
singularity at $x=0$, the integral of the left hand side of
(\ref{L*}) over $\Omega_{\varepsilon}=\Omega\bigcap
\{x:|x|>\varepsilon\}$ converges to the left hand side of
(\ref{q112}). Hence, in order to prove (\ref{cons}), it remains to
show that
\[
\int_{\Omega_{\varepsilon}}L^*Q dt dx \rightarrow
0,~~\varepsilon\rightarrow0.
\]
The integral above is equal to
\begin{equation} \label{last1}
\int_s^{\tau}\int_{|x|=\varepsilon}[-\frac{1}{2}\frac{\partial
Q}{\partial r}+\frac{\partial g}{\partial r}Q]d\sigma dt,
\end{equation}
where $d\sigma$ is the element of the surface area of the sphere
$|x|=\varepsilon$. The convergence of (\ref{last1}) to zero
follows immediately from (\ref{last2})\qed

\begin{lemma} \label{lnnn}
The family of processes $y^T(t)$, $T \geq 1$, is tight.
\end{lemma}
We shall prove this lemma below. First, however, we formulate the
main result of this section.
\begin{theorem} \label{nn1cc}
The distributions of the processes $y^T(t)$ converge as $T
\rightarrow \infty$, weakly in the space $C([0,1], \mathbb{R}^3)$,
to the distribution of the $3$-dimensional Markov  process with
continuous trajectories. The transition densities for the limiting
Markov process are given by (\ref{q11}) and (\ref{q12}).
\end{theorem}
\proof The convergence of the finite dimensional distributions of
$y^T(t)$ to those of the Markov process follows from (\ref{qdefe})
and Lemma~\ref{qpr}. Since the family $y^T(t)$ is tight, there is
a modification of the Markov process which has continuous
trajectories. \qed
\\
\\
\noindent {\it Proof of Lemma~\ref{lnnn}.} To prove tightness it
is enough to demonstrate that for each $\eta, \varepsilon > 0$
there are $0 < \delta < 1$ and $T_0 \geq 1$ such that for all $u
\in [0,1]$ we have
\begin{equation} \label{crite}
\mathrm{P}_{\beta, T} (\sup_{u \leq s \leq \min(t+\delta, 1)}
|y^T(s) - y^T(u)| > \varepsilon) \leq \delta \eta,~~T \geq T_0.
\end{equation}
Let $\eta, \varepsilon >0$ be fixed. Let $\mathcal{E}_\delta$ be
the event that a continuous function $x: [0,T] \rightarrow
{\mathbb{R}}^3$ satisfies
\[
\sup_{t \leq \delta T,} |x(t) - x(0)|/\sqrt{T}
> \varepsilon/8.
\]
Using arguments similar to those leading to (\ref{sup22}), we can
show that (\ref{crite}) follows from
\begin{equation} \label{sup22xx}
\sup_{x \in \mathbb{R}^d, |x| = \varepsilon \sqrt{T}/4}
\mathrm{E}^{x}_{0,T} ( \chi_{ \mathcal{E}_\delta} \exp( \int_0^{T}
\beta v (x(t)) d t) )\leq \delta \eta,~~T \geq T_0.
\end{equation}

 Let
\[ \tau = \min(\delta T, \inf\{t \geq 0: |x(t) - x(0)| = \varepsilon
\sqrt{T}/8\} ),
\]
The expectation in (\ref{sup22xx}) can be estimated as follows
\[
\mathrm{E}^{x}_{0,T} ( \chi_{ \mathcal{E}_\delta} \exp(\int_0^{T}
\beta v (x(t)) d t ))\leq \mathrm{E}^{x}_{0,T} ( \chi_{
\mathcal{E}_\delta} \mathrm{E}^{x(\tau)}_{0,T -\tau} \exp(
\int_0^{T-\tau} \beta v (x(t)) d t ))
\]
We claim that
\begin{equation} \label{lst}
\mathrm{E}^{x(\tau)}_{0,T -\tau} \exp( \int_0^{T-\tau} \beta v
(x(t)) d t ) \leq \sup_{x \in \mathbb{R}^d, |x| \geq \varepsilon
\sqrt{T}/8 } \mathrm{E}^{x}_{0,T} \exp( \int_0^{T} \beta v (x(t))
d t ) \leq c(\varepsilon)
\end{equation}
for some constant $c(\varepsilon)$ for all sufficiently large $T$.
It then remains to choose $\delta$ such that $\mathrm{E}^{x}_{0,T}
( \chi_{ \mathcal{E}_\delta} ) \leq \delta \eta /c(\varepsilon)$,
and the estimate  (\ref{sup22xx}) will follow. The second
inequality in (\ref{lst}) easily follows from part (2) of
Theorem~\ref{thz} and the fact that the probability of reaching
the support of $v$ before time $T$ by a Brownian path starting at
a distance $\varepsilon \sqrt{T}/8$ away from the origin is of
order $O(T^{-1/2})$ if $d  =3$. \qed

\end{document}